\documentclass[a4paper,11pt]{article}

\makeatletter
\gdef\@fpheader{   }
\gdef\@journal{jhep}

\RequirePackage{amsmath}
\RequirePackage{amssymb}
\RequirePackage{epsfig}
\RequirePackage{CJK}
\RequirePackage{graphicx}
\RequirePackage[numbers,sort&compress]{natbib}
\RequirePackage{color}
\RequirePackage[colorlinks=true
,urlcolor=blue
,anchorcolor=blue
,citecolor=blue
,filecolor=blue
,linkcolor=blue
,menucolor=blue
,pagecolor=blue
,linktocpage=true
,pdfproducer=medialab
,pdfa=true
,CJKbookmarks
]{hyperref}

\newif\ifnotoc\notocfalse
\newif\ifemailadd\emailaddfalse
\newif\iftoccontinuous\toccontinuousfalse

\def\@subheader{\@empty}
\def\@keywords{\@empty}
\def\@abstract{\@empty}
\def\@xtum{\@empty}
\def\@dedicated{\@empty}
\def\@arxivnumber{\@empty}
\def\@collaboration{\@empty}
\def\@collaborationImg{\@empty}
\def\@proceeding{\@empty}
\def\@preprint{\@empty}

\newcommand{\subheader}[1]{\gdef\@subheader{#1}}
\newcommand{\keywords}[1]{\if!\@keywords!\gdef\@keywords{#1}\else%
\PackageWarningNoLine{\jname}{Keywords already defined.\MessageBreak Ignoring last definition.}\fi}
\renewcommand{\abstract}[1]{\gdef\@abstract{#1}}
\newcommand{\dedicated}[1]{\gdef\@dedicated{#1}}
\newcommand{\arxivnumber}[1]{\gdef\@arxivnumber{#1}}
\newcommand{\proceeding}[1]{\gdef\@proceeding{#1}}
\newcommand{\xtumfont}[1]{\textsc{#1}}
\newcommand{\correctionref}[3]{\gdef\@xtum{\xtumfont{#1} \href{#2}{#3}}}
\newcommand\jname{JHEP}
\newcommand\acknowledgments{\section*{Acknowledgments}}

\newcommand\preprint[1]{\gdef\@preprint{\hfill #1}}

\newcommand\note[2][]{%
\if!#1!%
\stepcounter{footnote}\footnotetext{#2}%
\else%
{\renewcommand\thefootnote{#1}%
\footnotetext{#2}}%
\fi}

\newtoks\auth@toks
\renewcommand{\author}[2][]{%
  \if!#1!%
    \auth@toks=\expandafter{\the\auth@toks#2\ }%
  \else
    \auth@toks=\expandafter{\the\auth@toks#2$^{#1}$\ }%
  \fi
}

\newtoks\affil@toks\newif\ifaffil\affilfalse
\newcommand{\affiliation}[2][]{%
\affiltrue
  \if!#1!%
    \affil@toks=\expandafter{\the\affil@toks{\item[]#2}}%
  \else
    \affil@toks=\expandafter{\the\affil@toks{\item[$^{#1}$]#2}}%
  \fi
}

\newtoks\email@toks\newcounter{email@counter}%
\newcommand{\emailAdd}[1]{%
\emailaddtrue%
\ifnum\theemail@counter>0\email@toks=\expandafter{\the\email@toks, \@email{#1}}%
\else\email@toks=\expandafter{\the\email@toks\@email{#1}}%
\fi\stepcounter{email@counter}}
\newcommand{\@email}[1]{\href{mailto:#1}{\tt #1}}

\newcommand*\collaboration[1]{\gdef\@collaboration{#1}}
\newcommand*\collaborationImg[2][]{\gdef\@collaborationImg{#2}}

\newcommand\afterLogoSpace{\smallskip}
\newcommand\afterSubheaderSpace{\vskip3pt plus 2pt minus 1pt}
\newcommand\afterProceedingsSpace{\vskip21pt plus0.4fil minus15pt}
\newcommand\afterTitleSpace{\vskip23pt plus0.06fil minus13pt}
\newcommand\afterRuleSpace{\vskip23pt plus0.06fil minus13pt}
\newcommand\afterCollaborationSpace{\vskip3pt plus 2pt minus 1pt}
\newcommand\afterCollaborationImgSpace{\vskip3pt plus 2pt minus 1pt}
\newcommand\afterAuthorSpace{\vskip5pt plus4pt minus4pt}
\newcommand\afterAffiliationSpace{\vskip3pt plus3pt}
\newcommand\afterEmailSpace{\vskip16pt plus9pt minus10pt\filbreak}
\newcommand\afterXtumSpace{\par\bigskip}
\newcommand\afterAbstractSpace{\vskip16pt plus9pt minus13pt}
\newcommand\afterKeywordsSpace{\vskip16pt plus9pt minus13pt}
\newcommand\afterArxivSpace{\vskip3pt plus0.01fil minus10pt}
\newcommand\afterDedicatedSpace{\vskip0pt plus0.01fil}
\newcommand\afterTocSpace{\bigskip\medskip}
\newcommand\afterTocRuleSpace{\bigskip\bigskip}
\newlength{\affiliationsSep}
\newcommand\beforetochook{\pagestyle{myplain}\pagenumbering{roman}}

\DeclareFixedFont\trfont{OT1}{phv}{b}{sc}{11}

\renewcommand\maketitle{
\pagestyle{empty}
\thispagestyle{titlepage}
\setcounter{page}{0}
\noindent{\small\scshape\@fpheader}\@preprint\par
\afterLogoSpace
\if!\@subheader!\else\noindent{\trfont{\@subheader}}\fi
\afterSubheaderSpace
\if!\@proceeding!\else\noindent{\sc\@proceeding}\fi
\afterProceedingsSpace
{\LARGE\flushleft\sffamily\bfseries\@title\par}
\afterTitleSpace
\hrule height 1.5\p@%
\afterRuleSpace
\if!\@collaboration!\else
{\Large\bfseries\sffamily\raggedright\@collaboration}\par
\afterCollaborationSpace
\fi
\if!\@collaborationImg!\else
{\normalsize\bfseries\sffamily\raggedright\@collaborationImg}\par
\afterCollaborationImgSpace
\fi
{\bfseries\raggedright\sffamily\the\auth@toks\par}
\afterAuthorSpace
\ifaffil\begin{list}{}{%
\setlength{\leftmargin}{0.28cm}%
\setlength{\labelsep}{0pt}%
\setlength{\itemsep}{\affiliationsSep}%
\setlength{\topsep}{-\parskip}}
\itshape\small%
\the\affil@toks
\end{list}\fi
\afterAffiliationSpace
\ifemailadd 
\noindent\hspace{0.28cm}\begin{minipage}[l]{.9\textwidth}
\begin{flushleft}
\textit{E-mail:} \the\email@toks
\end{flushleft}
\end{minipage}
\else 
\PackageWarningNoLine{\jname}{E-mails are missing.\MessageBreak Plese use \protect\emailAdd\space macro to provide e-mails.}
\fi
\afterEmailSpace
\if!\@xtum!\else\noindent{\@xtum}\afterXtumSpace\fi
\if!\@abstract!\else\noindent{\renewcommand\baselinestretch{.9}\textsc{Abstract:}}\ \@abstract\afterAbstractSpace\fi
\if!\@keywords!\else\noindent{\textsc{Keywords:}} \@keywords\afterKeywordsSpace\fi
\if!\@arxivnumber!\else\noindent{\textsc{ArXiv ePrint:}} \href{http://arxiv.org/abs/\@arxivnumber}{\@arxivnumber}\afterArxivSpace\fi
\if!\@dedicated!\else\vbox{\small\it\raggedleft\@dedicated}\afterDedicatedSpace\fi
\ifnotoc\else
\iftoccontinuous\else\newpage\fi
\beforetochook\hrule
\tableofcontents
\afterTocSpace
\hrule
\afterTocRuleSpace
\fi
\setcounter{footnote}{0}
\pagestyle{myplain}\pagenumbering{arabic}
} 

\renewcommand{\baselinestretch}{1.1}\normalsize
\divide\@tempdima\baselineskip
\@tempcnta=\@tempdima
\skip\@mpfootins = \skip\footins
\renewcommand{\@dotsep}{10000}

\newcommand\ps@myplain{
\pagenumbering{arabic}
\renewcommand\@oddfoot{\hfill-- \thepage\ --\hfill}
\renewcommand\@oddhead{}}
\let\ps@plain=\ps@myplain

\newcommand\ps@titlepage{\renewcommand\@oddfoot{}\renewcommand\@oddhead{}}

\numberwithin{equation}{section}

\renewcommand\section{\@startsection{section}{1}{\z@}%
                                   {-3.5ex \@plus -1.3ex \@minus -.7ex}%
                                   {2.3ex \@plus.4ex \@minus .4ex}%
                                   {\normalfont\large\bfseries}}
\renewcommand\subsection{\@startsection{subsection}{2}{\z@}%
                                   {-2.3ex\@plus -1ex \@minus -.5ex}%
                                   {1.2ex \@plus .3ex \@minus .3ex}%
                                   {\normalfont\normalsize\bfseries}}
\renewcommand\subsubsection{\@startsection{subsubsection}{3}{\z@}%
                                   {-2.3ex\@plus -1ex \@minus -.5ex}%
                                   {1ex \@plus .2ex \@minus .2ex}%
                                   {\normalfont\normalsize\bfseries}}
\renewcommand\paragraph{\@startsection{paragraph}{4}{\z@}%
                                   {1.75ex \@plus1ex \@minus.2ex}%
                                   {-1em}%
                                   {\normalfont\normalsize\bfseries}}
\renewcommand\subparagraph{\@startsection{subparagraph}{5}{\parindent}%
                                   {1.75ex \@plus1ex \@minus .2ex}%
                                   {-1em}%
                                   {\normalfont\normalsize\bfseries}}

\def\fnum@figure{\textbf{\figurename\nobreakspace\thefigure}}
\def\fnum@table{\textbf{\tablename\nobreakspace\thetable}}

\long\def\@makecaption#1#2{%
  \vskip\abovecaptionskip
  \sbox\@tempboxa{\small #1. #2}%
  \ifdim \wd\@tempboxa >\hsize
    \small #1. #2\par
  \else
    \global \@minipagefalse
    \hb@xt@\hsize{\hfil\box\@tempboxa\hfil}%
  \fi
  \vskip\belowcaptionskip}

\renewenvironment{thebibliography}[1]{%
\begin{oldthebibliography}{#1}%
\small%
\raggedright%
\setlength{\itemsep}{5pt plus 0.2ex minus 0.05ex}%
}%
{%
\end{oldthebibliography}%
}

\makeatother

\usepackage[T1]{fontenc} 
\usepackage{romannum} 
\usepackage{CJK}

\title{{\boldmath Renormalization of divergent moment in probability theory}}

\author[a,1]{Ping Zhang,}\note{zhangping@cueb.edu.cn.}
\author[b,2]{Wen-Du Li,}\note{liwendu@tjnu.edu.cn.}
\author[c,3]{and Wu-Sheng Dai}\note{daiwusheng@tju.edu.cn.}

\affiliation[a]{School of Finance, Capital University of Economics and Business, Beijing 100070, PR China}
\affiliation[b]{College of Physics and Materials Science, Tianjin Normal University, Tianjin 300387, PR China}
\affiliation[c]{Department of Physics, Tianjin University, Tianjin 300350, P.R. China}

\abstract{
Some probability distributions have moments, and some do not. For example, the
normal distribution has power moments of arbitrary order, but the Cauchy
distribution does not have power moments. In this paper, by analogy with the
renormalization method in quantum field theory, we suggest a renormalization
scheme to remove the divergence in divergent moments. We establish more than
one renormalization procedure to renormalize the same moment to prove that the
renormalized moment is scheme-independent. The power moment is usually a
positive-integer-power moment; in this paper, we introduce
nonpositive-integer-power moments by a similar treatment of renormalization.
An approach to calculating logarithmic moment from power moment is proposed,
which can serve as a verification of the validity of the renormalization
procedure. The renormalization schemes proposed are the zeta function scheme,
the subtraction scheme, the weighted moment scheme, the cut-off scheme, the
characteristic function scheme, the Mellin transformation scheme, and the
power-logarithmic moment scheme. The probability distributions considered are
the Cauchy distribution, the Levy distribution, the $q$-exponential
distribution, the $q$-Gaussian distribution, the normal distribution, the
Student's $t$-distribution, and the Laplace distribution. 

}

\begin{document} 
\begin{CJK*}{GBK}{song}
\maketitle 

\flushbottom

\section{Introduction}

The moment of a distribution with the probability density function $p\left(
x\right)  $, generally, is defined by%
\begin{equation}
m_{f}=\int_{-\infty}^{\infty}p\left(  x\right)  f\left(  x\right)  dx.
\label{moment}%
\end{equation}
Different $f\left(  x\right)  $ define different moments. The most familiar
moments are the $n$-th power moment and the logarithmic moment, corresponding
to $f\left(  x\right)  $ is a power function or a logarithmic function.

Choosing $f\left(  x\right)  =x^{n}$ defines the $n$-th power moment:%
\begin{equation}
m_{n}=\int_{-\infty}^{\infty}p\left(  x\right)  x^{n}dx. \label{norder}%
\end{equation}
The zeroth moment corresponds to the normalization, the first moment is the
mean value, the second central moment is the variance, the third central
moment is the skewness, and the fourth central moment is the kurtosis
\cite{rice2006mathematical}. Choosing $f\left(  x\right)  =\ln^{n}x$ defines
the $n$-th logarithmic moment \cite{nicolas2002introduction}:%
\begin{equation}
\widetilde{m}_{n}=\int_{-\infty}^{\infty}p\left(  x\right)  \ln^{n}xdx.
\label{logM}%
\end{equation}
The moment of continuous distributions is usually defined by integrals. If the
function $p\left(  x\right)  f\left(  x\right)  $ in Eq. (\ref{moment}) is
nonintegrable, the moment $m_{f}$ does not exist. The existence of the $n$-th
moment, by definition (\ref{norder}), relies on whether $p\left(  x\right)
x^{n}$ is integrable or not. Similarly, the existence of the logarithmic
moment requires $p\left(  x\right)  \ln^{n}x$ to be integrable. In various
statistical distributions, some have moments, such as the normal distribution,
the Laplace distribution, and Student's $t$-distribution; some have no
moments, such as the Cauchy distribution and the Levy distribution; some have
moments only under some special values of parameters, such as the
$q$-exponential distribution and the $q$-Gaussian distribution
\cite{johnson1995continuous}. In this paper, by analogy with the
renormalization method in quantum field theory, we suggest a renormalization
scheme to remove the divergence in the integral definition of moments, Eq.
(\ref{moment}), to achieve a finite renormalized moment.

For the distribution that has no moment, the integral in the definition
(\ref{moment}) is divergent. The renormalization treatment aims to achieve a
finite moment by removing the divergence in the integral definition
(\ref{moment}). Roughly speaking, such a treatment is to obtain a finite value
by subtracting infinity from infinity. However, this treatment will naturally
cause a problem: whether the renormalized result depends on the choice of
subtraction scheme. After all, infinity minus infinity can give any value. A
correct subtraction treatment must be subtraction-scheme independent. The
renormalized results given by the different renormalization schemes must be
the same.

In order to show the validity of the renormalization schemes suggested in this
paper, we construct a variety of renormalization schemes and renormalize an
$n$-th moment with different renormalization schemes. If the renormalized
$n$-th moment given by different schemes are the same, then the
renormalization-scheme independence is verified.

Constructing more than one renormalization scheme also has a technical reason.
The renormalization procedure, practically, depends on exact results, such as
exact integrals and exact sums. Nevertheless, the exact result is not always
available. We will have more opportunities if we have more than one
renormalization scheme. The exact result that cannot be obtained by one scheme
may be obtained by other schemes. Therefore, more renormalization schemes mean
more operability.

In this paper, we will construct the following renormalization schemes: zeta
function scheme, subtraction scheme, weighted moment scheme, cut-off scheme,
characteristic function scheme, Mellin transformation scheme, and
power-logarithmic moment scheme.

Using the renormalization scheme suggested in this paper, we will calculate
the renormalized $n$-th power moment for the following distributions: Cauchy
distribution, Levy distribution, $q$-exponential distribution, and
$q$-Gaussian distribution. We use more than one renormalization scheme for
each distribution. Whether a renormalization scheme applies to a certain
distribution depends on, for example, whether the integral in the
renormalization procedure can be exactly worked out.

For power moments, we usually only consider positive-integer power moments.
However, in the renormalization treatment, we can obtain arbitrary real-number
power moments and even complex-number power moments. By observing a
complex-number power moment, we can analyze the singularity structure of the
moment function on the complex plane. For example, the Laplace distribution
has any positive-integer power moments. However, when we extend the moment
function of the Laplace distribution to the complex plane, we will find that
the complex moment function of the Laplace distribution is singular at the
negative-integer power, and the negative-integer power is the singularity of
complex-number moment function. The Laplace distribution does not have
negative-integer power moments. In order to obtain the negative-integer power
moments of the Laplace distribution, we need to use the renormalization method
suggested in this paper.

In addition to renormalization schemes, we also provide a method to calculate
logarithmic moments from power moments. The method is suitable for the case
that the calculation of logarithmic moments is more complicated than that of
power moments. However, some distributions have logarithmic moments but have
no power moments. That is, the integral in the logarithmic moment (\ref{logM})
converges, but the integral in the power moment (\ref{norder}) diverges. If we
still want to use the power moment to calculate the logarithmic moment, we can
use the renormalized power moment instead of the divergent power moment. In
this paper, we give some examples of calculating logarithmic moments from
power moments, including Cauchy distribution, Levy distribution,
$q$-exponential distribution, $q$-Gaussian distribution, normal distribution,
Student's $t$-distribution, and Laplace distribution. In these examples, some
distributions have power moments, and others do not. For distributions that
have no power moments, the renormalized power moment needs to be calculated first.

Calculating logarithmic moments from power moments also verifies the validity
of the renormalization scheme. For a distribution that has no power moments
but has logarithmic moments, we only need to renormalize the power moment, and
the logarithmic moment can be obtained directly by definition (\ref{logM}).
Comparing the logarithmic moment obtained from the renormalized power moment
with the logarithmic moment calculated by definition (\ref{logM}), we can
verify the validity of renormalization.

The Cauchy distribution is widely used in physics and mathematics, such as the
phase oscillator systems \cite{tanaka2020low}, the fractional generalized
Cauchy process \cite{uchiyama2019fractional}, many-body localized systems
\cite{filippone2016drude}, the edge of chaos and avalanches in neural networks
\cite{kusmierz2020edge},\ the factor models and contingency tables
\cite{pillai2016unexpected}, the Dirichlet random probability
\cite{letac2018dirichlet}, the statistical depth \cite{einmahl2015bridging},
and the image processing \cite{mei2018cauchy}. The Levy distribution is a
special case of the Levy $\alpha$-stable distribution \cite{rocha2019levy}.
The Levy distribution plays important roles in the geometry of multiparameter
families of quantum states \cite{penner2021hilbert}, the efficiency of random
target search strategies quantification \cite{levernier2021universality}, the
nondiffusive suprathermal ion transport \cite{manke2019truncated}, and the
multimode viscous hydrodynamics for one-dimensional spinless electrons
\cite{protopopov2021anomalous}. The $q$-exponential distribution and the
$q$-Gaussian distribution are examples of the Tsallis distributions arising
from the maximization of the Tsallis entropy under appropriate constraints
\cite{tsallis2009introduction}. The $q$-exponential distribution is a
generalization of the exponential distribution, and the $q$-Gaussian
distribution is a generalization of the Gaussian distribution
\cite{domingo2017boundedness}.\ The $q$-exponential distribution and the
$q$-Gaussian distribution are popular in complex systems
\cite{mehri2011keyword}, the two-qubit quantum system and the
transverse-momentum behavior of hadrons in high-energy proton-proton
collisions \cite{jizba2019maximum}, the quantum Coulomb system in a confining
potential \cite{bonart2020scaling}, financial market
\cite{alonso2019q,liu2021european,liu2012intermediate,denys2016universality},
the clustering algorithm for cryo-electron microscopy \cite{chen2014gamma},
the Wasserstein Geometry \cite{takatsu2013behaviors}, and the Swarm
Quantum-like Particle Optimization algorithms \cite{kamberaj2014q}.

In section \ref{powerM}, we construct renormalization schemes for power
moments. In section \ref{RenEx}, we apply the renormalization scheme to
calculate the renormalized power moment for various distributions. In section
\ref{noninteger}, we discuss the generalization of the renormalized power
moment. We provide a method to calculate logarithmic moments from power
moments in section \ref{mLog}, and give examples in section \ref{mLogEx}.
Discussions and conclusion are in section \ref{conclusion}.

\section{Renormalization scheme: $n$-th power moment \label{powerM}}

In this section, to remove divergences in power moments of distributions that
have no well-defined moment, we construct several different renormalization
schemes, including the zeta function scheme, the subtraction scheme, the
weighted moment scheme, the cut-off scheme, the characteristic function
scheme, the Mellin transformation scheme, and the power-logarithmic moment
scheme. Moreover, in a renormalization scheme, we introduce the
power-logarithmic moment. In the following, we apply these renormalization
schemes to various distributions.

\subsection{Zeta function scheme \label{ZetaR}}

The key step in the renormalization procedure is the analytical continuation.
The analytical continuation of the zeta function is well studied in
mathematics and has been systematically applied to the renormalization in
quantum field theory \cite{vassilevich2003heat}. In this section, we establish
the zeta function renormalization scheme for the distribution that has no
power moment.

A spectrum is a series of numbers, discrete or continuous. For example, all
eigenvalues of an operator form an eigenvalue spectrum. From a given spectrum,
we can define spectral functions, such as spectral zeta functions, one-loop
effective actions, vacuum energies, and various thermodynamic quantities
\cite{dai2009number,dai2010approach}. Here, we focus on the spectral zeta function.

\textit{Generalized} \textit{spectral zeta function. }For a discrete spectrum
$\left\{  \lambda_{0},\lambda_{1},\ldots,\lambda_{n},\ldots\right\}  $, the
spectral zeta function is defined as $\zeta\left(  s\right)  =\sum
_{\lambda_{n}\geq\lambda_{0}\left(  \neq-\infty\right)  }\lambda_{n}^{-s}$
\cite{dai2010approach}. For a continuous spectrum, with the spectral density
function $\rho\left(  \lambda\right)  $, the spectral zeta function is defined
as $\zeta\left(  s\right)  =\int_{\lambda_{0}}^{\infty}d\lambda\rho\left(
\lambda\right)  \lambda^{-s}$. In physics, however, the physical operators are
usually lower bounded. The above definition of spectral zeta functions applies
only to lower bounded spectra. Since the random variables in many
distributions are ranged from $-\infty$ to $+\infty$, before we apply the
spectral zeta function to statistical distributions, we must generalize the
spectral zeta function to the spectrum that is neither lower bounded nor upper bounded.

For the spectrum that is neither lower bounded nor upper bounded, for the
discrete spectrum, we generalize the spectral zeta function as%
\begin{equation}
\zeta\left(  s\right)  =\sum_{-\infty<\lambda_{n}<\infty}\lambda_{n}^{-s}.
\label{zetadef}%
\end{equation}
and for the continuous spectrum,
\begin{equation}
\zeta\left(  s\right)  =\int_{-\infty}^{\infty}d\lambda\rho\left(
\lambda\right)  \lambda^{-s}. \label{zetaint}%
\end{equation}

\textit{Moment as zeta function.} To perform renormalization with the zeta
function, we regard that all random variables form a spectrum. In this way,
each probability distribution has its own zeta function.

Concretely, we regard the random variable $x$ in the distribution function as
an element $\lambda$ in a spectrum, and regard the probability density
function $p\left(  x\right)  $ as the spectral density function $\rho\left(
\lambda\right)  $. Thus, taking continuous spectra as an example, the spectral
zeta function (\ref{zetaint}) becomes%

\begin{equation}
\zeta\left(  s\right)  =\int_{-\infty}^{\infty}dxp\left(  x\right)  x^{-s}.
\end{equation}
Comparing with the definition of the $n$-th power moment (\ref{norder}) shows
that the $n$-th power moment is a spectral zeta function:
\begin{equation}
m_{n}^{\text{ren}}=\zeta\left(  -n\right)  , \label{mnzeta}%
\end{equation}
where $m_{n}^{\text{ren}}$ denotes the $n$-th renormalized power moment.

It should be emphasized that although the integral in the definition of the
$n$-th moment (\ref{norder}) and the spectral zeta function (\ref{zetaint})
are the same in form, the analytic region of the spectral zeta function is
often larger than the integrable region of the $n$-th moment. This is because
different representations of the zeta function have different analytic
regions. When we express the $n$-th moment with the spectral zeta function, we
have indeed made an analytic continuation and extended the analytic region. In
this way, the original divergent moment becomes convergent, and the divergent
moment no longer diverges.

Although analytic continuation can remove the divergence caused by the
integral definition of the moment, the true singularity in the moment
function, the obstacle of analytic continuation, can not be removed by
analytic continuation. For true singularities, the role of analytical
continuation is to expose singularities. These true singularities will be
removed by the minimal subtraction discussed later.

\subsection{Subtraction method \label{remove}}

Consider a divergent integral%
\begin{equation}
m=\int_{-\infty}^{\infty}f\left(  x\right)  dx. \label{mInt}%
\end{equation}
We remove divergence by introducing counterterm terms. We consider the
divergence of the integral at $x=0$, $x\rightarrow\infty$, $x\rightarrow
-\infty$, and $x\rightarrow\pm\infty$, respectively.

\subsubsection{Divergence at $x=0$ \label{remove0}}

First, consider the case that the integral diverges at the lower limit of the
integral, $x=0$.

Expand the integrand at $x=0$:%
\begin{align}
f\left(  x\right)   &  =\sum_{j=1}^{\infty}a_{j}x^{\alpha_{j}}\nonumber\\
&  =\sum_{j=1}^{j_{\max}}a_{j}x^{\alpha_{j}}+\sum_{j_{\min}}^{\infty}%
a_{j}x^{\alpha_{j}}, \label{series}%
\end{align}
where $j_{\max}$ satisfies $\operatorname{Re}\alpha_{j}\leq-1$ and $j_{\min}$
satisfies $\operatorname{Re}\alpha_{j}>-1$. The divergence at the lower limit
$x=0$ comes from $\operatorname{Re}\alpha_{j}\leq-1$. Thus, we use the method
in Ref. \cite{zeidler2008quantum} to remove the divergence.

Rewrite the integral (\ref{mInt}) as
\begin{align}
m  &  =\left[  \int_{-1}^{0}f\left(  x\right)  dx+\int_{0}^{1}f\left(
x\right)  dx\right]  +\left[  \int_{-\infty}^{-1}f\left(  x\right)
dx+\int_{1}^{\infty}f\left(  x\right)  dx\right] \nonumber\\
&  =I_{\text{div}}+I_{\text{con}}. \label{mdc}%
\end{align}
The divergence occurs in $I_{\text{div}}=\int_{-1}^{0}f\left(  x\right)
dx+\int_{0}^{1}f\left(  x\right)  dx$, and the divergence of $I_{\text{div}}$
comes from the $\operatorname{Re}\alpha_{j}\leq-1$ term $\sum_{j=1}^{j_{\max}%
}a_{j}x^{\alpha_{j}}$. Next, we subtract out the divergence term in
$I_{\text{div}}$.

Write $I_{\text{div}}$\ as
\begin{align}
I_{\text{div}}  &  =\int_{-1}^{0}\left[  f\left(  x\right)  -\sum
_{j=1}^{j_{\max}}a_{j}x^{\alpha_{j}}\right]  dx+\int_{0}^{1}\left[  f\left(
x\right)  -\sum_{j=1}^{j_{\max}}a_{j}x^{\alpha_{j}}\right]  dx\nonumber\\
&  +\int_{-1}^{0}\left[  \sum_{j=1}^{j_{\max}}a_{j}x^{\alpha_{j}}\right]
dx+\int_{0}^{1}\left[  \sum_{j=1}^{j_{\max}}a_{j}x^{\alpha_{j}}\right]  dx.
\end{align}
Here, $I_{\text{div}}$\ consists of four integrals. The first two integrals
are already convergent because the divergent part ($\operatorname{Re}a_{j}%
\leq-1$) has been subtracted out. The divergence of $I_{\text{div}}$ is in the
last two integrals.

In the divergent integral, a negative-power term $a_{-1}x^{-1}$ may appear and
the integral of this term is a logarithmic function. We write this term
separately as follows:
\begin{align}
&  \int_{-1}^{0}\left[  \sum_{j=1}^{j_{\max}}a_{j}x^{\alpha_{j}}\right]
dx\nonumber\\
&  =\left[  \underset{j=1}{\overset{j_{\max}}{\sum\nolimits^{\prime}}}%
\frac{a_{j}}{\alpha_{j}+1}x^{\alpha_{j}+1}+a_{-1}\ln x\right]  _{x=0^{-}%
}-\left[  \underset{j=1}{\overset{j_{\max}}{\sum\nolimits^{\prime}}}%
\frac{a_{j}}{\alpha_{j}+1}x^{\alpha_{j}+1}+a_{-1}\ln x\right]  _{x=-1},
\end{align}
and
\begin{align}
&  \int_{0}^{1}\left[  \sum_{j=1}^{j_{\max}}a_{j}x^{\alpha_{j}}\right]
dx\nonumber\\
&  =\left[  \underset{j=1}{\overset{j_{\max}}{\sum\nolimits^{\prime}}}%
\frac{a_{j}}{\alpha_{j}+1}x^{\alpha_{j}+1}+a_{-1}\ln x\right]  _{x=1}-\left[
\underset{j=1}{\overset{j_{\max}}{\sum\nolimits^{\prime}}}\frac{a_{j}}%
{\alpha_{j}+1}x^{\alpha_{j}+1}+a_{-1}\ln x\right]  _{x=0^{+}},
\end{align}
where $\sum\nolimits^{^{\prime}}$\ denotes that the negative-power term is not
included in the sum, i.e., $\alpha_{j}\neq-1$. Obviously, these two integrals
converge at $x=-1$ and $x=1$, and the divergence comes from $x=0$. After
removing the divergent term $\left[  \underset{j=1}{\overset{j_{\max}%
}{\sum\nolimits^{\prime}}}\frac{a_{j}}{\alpha_{j}+1}x^{\alpha_{j}+1}+a_{-1}\ln
x\right]  _{x=0^{\pm}}$, we arrive at a renormalized result,
\begin{align}
I_{\text{div}}^{\text{ren}}  &  =\int_{-1}^{0}\left[  f\left(  x\right)
-\sum_{j=1}^{j_{\max}}a_{j}x^{\alpha_{j}}\right]  dx+\int_{0}^{1}\left[
f\left(  x\right)  -\sum_{j=1}^{j_{\max}}a_{j}x^{\alpha_{j}}\right]
dx\nonumber\\
&  -\underset{j=1}{\overset{j_{\max}}{\sum\nolimits^{\prime}}}\frac{\left(
-1\right)  ^{\alpha_{j}+1}a_{j}}{\alpha_{j}+1}+i\pi a_{-1}%
+\underset{j=1}{\overset{j_{\max}}{\sum\nolimits^{\prime}}}\frac{a_{j}}%
{\alpha_{j}+1}. \label{Idren}%
\end{align}

The renormalized $n$-th moment, by Eqs. (\ref{mdc}) and (\ref{Idren}), reads
\begin{align}
m^{\text{ren}}  &  =\int_{-\infty}^{-1}f\left(  x\right)  dx+\int_{-1}%
^{0}\left[  f\left(  x\right)  -\sum_{j=1}^{j_{\max}}a_{j}x^{\alpha_{j}%
}\right]  dx+\int_{0}^{1}\left[  f\left(  x\right)  -\sum_{j=1}^{j_{\max}%
}a_{j}x^{\alpha_{j}}\right]  dx\nonumber\\
&  -\underset{j=1}{\overset{j_{\max}}{\sum\nolimits^{\prime}}}\frac{\left(
-1\right)  ^{\alpha_{j}+1}a_{j}}{\alpha_{j}+1}+i\pi a_{-1}%
+\underset{j=1}{\overset{j_{\max}}{\sum\nolimits^{\prime}}}\frac{a_{j}}%
{\alpha_{j}+1}+\int_{1}^{\infty}f\left(  x\right)  dx.
\end{align}

\subsubsection{Divergence at $x\rightarrow\infty$ \label{removeinf}}

If the integral diverges at $x\rightarrow\infty$, we expand the integrand at
$x\rightarrow\infty$:%
\begin{equation}
f\left(  x\right)  =\sum_{j=1}^{\infty}b_{j}\frac{1}{x^{\beta_{j}}}.
\end{equation}
The terms of $\operatorname{Re}\beta_{j}\leq1$\ will diverge. We use the
treatment above to remove the divergence.

Write the integral (\ref{mInt}) in two parts:
\begin{align}
m  &  =\int_{-\infty}^{1}f\left(  x\right)  dx+\int_{1}^{\infty}f\left(
x\right)  dx\nonumber\\
&  =I_{\text{con}}+I_{\text{div}}. \label{minf}%
\end{align}
$I_{\text{con}}$ is convergent and $I_{\text{div}}$ diverges.

Write the divergent integral $I_{\text{div}}$\ as%
\begin{align}
I_{\text{div}}  &  =\int_{1}^{\infty}\left[  f\left(  x\right)  -\sum
_{j=1}^{j_{\max}}b_{j}\frac{1}{x^{\beta_{j}}}\right]  dx+\int_{1}^{\infty
}\left[  \sum_{j=1}^{j_{\max}}b_{j}\frac{1}{x^{\beta_{j}}}\right]
dx\nonumber\\
&  =\int_{1}^{\infty}\left[  f\left(  x\right)  -\sum_{j=1}^{j_{\max}}%
b_{j}\frac{1}{x^{\beta_{j}}}\right]  dx+\left[
\underset{j=1}{\overset{j_{\max}}{\sum\nolimits^{\prime}}}\frac{b_{j}}%
{1-\beta_{j}}\frac{1}{x^{\beta_{j}-1}}+b_{-1}\ln x\right]  _{x\rightarrow
\infty}\nonumber\\
&  -\left[  \underset{j=1}{\overset{j_{\max}}{\sum\nolimits^{\prime}}}%
\frac{b_{j}}{1-\beta_{j}}\frac{1}{x^{\beta_{j}-1}}+b_{-1}\ln x\right]  _{x=1},
\end{align}
where the sum $\sum\nolimits^{^{\prime}}$ does not contain the negative-power
term, i.e., $\beta_{j}\neq1$, and $j_{\max}$ satisfies $\operatorname{Re}%
\beta_{j}\leq1$. After removing the divergent term $\left[
\underset{j=1}{\overset{j_{\max}}{\sum\nolimits^{\prime}}}\frac{b_{j}}%
{1-\beta_{j}}\frac{1}{x^{\beta_{j}-1}}+b_{-1}\ln x\right]  _{x\rightarrow
\infty}$, we arrive at a renormalized result
\begin{equation}
I_{\text{div}}^{\text{ren}}=\int_{1}^{\infty}\left[  f\left(  x\right)
-\sum_{j=1}^{j_{\max}}b_{j}\frac{1}{x^{\beta_{j}}}\right]
dx-\underset{j=1}{\overset{j_{\max}}{\sum\nolimits^{\prime}}}\frac{b_{j}%
}{1-\beta_{j}}. \label{Iinfdivren}%
\end{equation}

The renormalized $n$-th moment, by Eqs. (\ref{minf}) and (\ref{Iinfdivren}),
reads%
\begin{equation}
m^{\text{ren}}=\int_{-\infty}^{1}f\left(  x\right)  dx+\int_{1}^{\infty
}\left[  f\left(  x\right)  -\sum_{j=1}^{j_{\max}}b_{j}\frac{1}{x^{\beta_{j}}%
}\right]  dx-\underset{j=1}{\overset{j_{\max}}{\sum\nolimits^{\prime}}}%
\frac{b_{j}}{1-\beta_{j}}. \label{mreninf}%
\end{equation}

\subsubsection{Divergence at $x\rightarrow-\infty$ \label{removeinf-}}

If the integral diverges at $x\rightarrow-\infty$, we expand the integrand at
$x\rightarrow-\infty$:%
\begin{equation}
f\left(  x\right)  =\sum_{j=1}^{\infty}c_{j}\frac{1}{x^{\gamma_{j}}}.
\end{equation}
The terms of $\operatorname{Re}\gamma_{j}\leq1$\ will diverge.

Write the integral (\ref{mInt}) in two parts:%
\begin{align}
m  &  =\int_{-\infty}^{-1}f\left(  x\right)  dx+\int_{-1}^{\infty}f\left(
x\right)  dx\label{minf-}\\
&  =I_{\text{div}}+I_{\text{con}}.
\end{align}
$I_{\text{con}}$ is convergent and $I_{\text{div}}$ diverges.

Write the divergent integral $I_{\text{div}}$\ as%
\begin{align}
I_{\text{div}}  &  =\int_{-\infty}^{-1}\left[  f\left(  x\right)  -\sum
_{j=1}^{j_{\max}}c_{j}\frac{1}{x^{\gamma_{j}}}\right]  dx+\int_{-\infty}%
^{-1}\left[  \sum_{j=1}^{j_{\max}}c_{j}\frac{1}{x^{\gamma_{j}}}\right]
dx\nonumber\\
&  =\int_{-\infty}^{-1}\left[  f\left(  x\right)  -\sum_{j=1}^{j_{\max}}%
c_{j}\frac{1}{x^{\gamma_{j}}}\right]  dx+\left[
\underset{j=1}{\overset{j_{\max}}{\sum\nolimits^{\prime}}}\frac{c_{j}%
}{1-\gamma_{j}}\frac{1}{x^{\gamma_{j}-1}}+c_{-1}\ln x\right]  _{x\rightarrow
-1}\nonumber\\
&  -\left[  \underset{j=1}{\overset{j_{\max}}{\sum\nolimits^{\prime}}}%
\frac{c_{j}}{1-\gamma_{j}}\frac{1}{x^{\gamma_{j}-1}}+c_{-1}\ln x\right]
_{x=-\infty},
\end{align}
where the sum $\sum\nolimits^{^{\prime}}$ does not contain the negative-power
term, i.e., $\gamma_{j}\neq1$, and $j_{\max}$ satisfies $\operatorname{Re}%
\gamma_{j}\leq1$. After removing the divergent term $\left[
\underset{j=1}{\overset{j_{\max}}{\sum\nolimits^{\prime}}}\frac{c_{j}%
}{1-\gamma_{j}}\left(  \frac{1}{x}\right)  ^{\gamma_{j}-1}+c_{-1}\ln x\right]
_{x=-\infty}$, we arrive at a renormalized result%
\begin{equation}
I_{\text{div}}^{\text{ren}}=\int_{-\infty}^{-1}\left[  f\left(  x\right)
-\sum_{j=1}^{j_{\max}}c_{j}\frac{1}{x^{\gamma_{j}}}\right]
dx+\underset{j=1}{\overset{j_{\max}}{\sum\nolimits^{\prime}}}\frac{\left(
-1\right)  ^{\gamma_{j}-1}c_{j}}{1-\gamma_{j}}+i\pi c_{-1}.
\label{Iinfdivren-}%
\end{equation}

The renormalized $n$-th moment, by Eqs. (\ref{minf}) and (\ref{Iinfdivren}),
reads%
\begin{equation}
m^{\text{ren}}=\int_{-1}^{\infty}f\left(  x\right)  dx+\int_{-\infty}%
^{-1}\left[  f\left(  x\right)  -\sum_{j=1}^{j_{\max}}c_{j}\frac{1}%
{x^{\gamma_{j}}}\right]  dx+\underset{j=1}{\overset{j_{\max}}{\sum
\nolimits^{\prime}}}\frac{\left(  -1\right)  ^{\gamma_{j}-1}c_{j}}%
{1-\gamma_{j}}+i\pi c_{-1}. \label{mreninf-}%
\end{equation}

\subsubsection{Divergence at $x\rightarrow\pm\infty$}

\bigskip If the integral diverges at both $x\rightarrow\pm\infty$, we expand
the integrand at $x\rightarrow\pm\infty$, respectively:%
\begin{equation}
f\left(  x\right)  =\sum_{j=1}^{\infty}d_{j}\frac{1}{x^{\delta_{j}}}.
\end{equation}
The terms of $\operatorname{Re}\delta_{j}\leq1$\ will diverge.

Write the integral (\ref{mInt}) in two parts:%
\begin{align}
m  &  =\int_{-1}^{1}f\left(  x\right)  dx+\left[  \int_{-\infty}^{-1}f\left(
x\right)  dx+\int_{1}^{\infty}f\left(  x\right)  dx\right] \nonumber\\
&  =I_{\text{con}}+I_{\text{div}}. \label{minf+-}%
\end{align}
$I_{\text{con}}$ is convergent and $I_{\text{div}}$ diverges.

Write the divergent integral $I_{\text{div}}$\ as%
\begin{align}
I_{\text{div}}  &  =\int_{-\infty}^{-1}\left[  f\left(  x\right)  -\sum
_{j=1}^{j_{\max}}d_{j}\frac{1}{x^{\delta_{j}}}\right]  dx+\int_{-\infty}%
^{-1}\left[  \sum_{j=1}^{j_{\max}}d_{j}\frac{1}{x^{\delta_{j}}}\right]
dx\nonumber\\
&  +\int_{1}^{\infty}\left[  f\left(  x\right)  -\sum_{j=1}^{j_{\max}}%
d_{j}\frac{1}{x^{\delta_{j}}}\right]  dx+\int_{1}^{\infty}\left[  \sum
_{j=1}^{j_{\max}}d_{j}\frac{1}{x^{\delta_{j}}}\right]  dx,
\end{align}
where $j_{\max}$ satisfies $\operatorname{Re}\delta_{j}\leq1$. Working out the
integral gives
\begin{align}
I_{\text{div}}  &  =\int_{-\infty}^{-1}\left[  f\left(  x\right)  -\sum
_{j=1}^{j_{\max}}d_{j}\frac{1}{x^{\delta_{j}}}\right]  dx+\left.
\underset{j=1}{\overset{j_{\max}}{\sum\nolimits^{\prime}}}\frac{d_{j}%
}{1-\delta_{j}}x^{1-\delta_{j}}\right\vert _{x=-\infty}^{x=-1}+\left.
d_{-1}\ln x\right\vert _{x=-\infty}^{x=-1}\nonumber\\
&  +\int_{1}^{\infty}\left[  f\left(  x\right)  -\sum_{j=1}^{j_{\max}}%
d_{j}\frac{1}{x^{\delta_{j}}}\right]  dx+\left.
\underset{j=1}{\overset{j_{\max}}{\sum\nolimits^{\prime}}}\frac{d_{j}%
}{1-\delta_{j}}x^{1-\delta_{j}}\right\vert _{x=1}^{x=\infty}+\left.  d_{-1}\ln
x\right\vert _{x=1}^{x=\infty},
\end{align}
where the sum $\sum\nolimits^{^{\prime}}$ does not contain the negative-power
term, i.e., $\delta_{j}\neq1$. After dropping the divergent terms, $\left.
\underset{j=1}{\overset{j_{\max}}{\sum\nolimits^{\prime}}}x^{1-\delta_{j}%
}\right\vert _{x=-\infty}$, $\left.  d_{-1}\ln x\right\vert _{x=-\infty}$,
$\left.  \underset{j=1}{\overset{j_{\max}}{\sum\nolimits^{\prime}}}\frac
{d_{j}}{1-\delta_{j}}x^{1-\delta_{j}}\right\vert _{x=\infty}$, and $\left.
d_{-1}\ln x\right\vert _{x=\infty}$, we arrive at a renormalized result
\begin{align}
I_{\text{div}}^{\text{ren}}  &  =\int_{-\infty}^{-1}\left[  f\left(  x\right)
-\sum_{j=1}^{j_{\max}}d_{j}\frac{1}{x^{\delta_{j}}}\right]
dx+\underset{j=1}{\overset{j_{\max}}{\sum\nolimits^{\prime}}}\frac{\left(
-1\right)  ^{1-\delta_{j}}d_{j}}{1-\delta_{j}}+id_{-1}\pi\nonumber\\
&  +\int_{1}^{\infty}\left[  f\left(  x\right)  -\sum_{j=1}^{j_{\max}}%
d_{j}\frac{1}{x^{\delta_{j}}}\right]  dx-\underset{j=1}{\overset{j_{\max
}}{\sum\nolimits^{\prime}}}\frac{d_{j}}{1-\delta_{j}}. \label{Iinfdivren+-}%
\end{align}

The renormalized $n$-th moment, by Eqs. (\ref{minf+-}) and (\ref{Iinfdivren+-}%
), reads%
\begin{align}
m^{\text{ren}}  &  =\int_{-1}^{1}f\left(  x\right)  dx+\int_{-\infty}%
^{-1}\left[  f\left(  x\right)  -\sum_{j=1}^{j_{\max}}d_{j}\frac{1}%
{x^{\delta_{j}}}\right]  dx+\underset{j=1}{\overset{j_{\max}}{\sum
\nolimits^{\prime}}}\frac{\left(  -1\right)  ^{1-\delta_{j}}d_{j}}%
{1-\delta_{j}}+id_{-1}\pi\nonumber\\
&  +\int_{1}^{\infty}\left[  f\left(  x\right)  -\sum_{j=1}^{j_{\max}}%
d_{j}\frac{1}{x^{\delta_{j}}}\right]  dx-\underset{j=1}{\overset{j_{\max
}}{\sum\nolimits^{\prime}}}\frac{d_{j}}{1-\delta_{j}}. \label{mreninf+-}%
\end{align}

\subsection{Weight function scheme \label{WFR}}

In this section, we establish a weighted moment scheme to renormalize power
moments. In the weighted moment scheme, we introduce a weighted moment and
then remove the divergence in the weighted moment by a renormalization
procedure. The weighted moment, of course, depends on the choice of weight
functions \cite{liu2012intermediate}, but the renormalized moment must be
independent of the weight function, i.e., the weighted moment scheme must be scheme-independent.

For the divergent $n$-th moment (\ref{norder}), we choose the weight function
$g_{s}\left(  x\right)  $ satisfying%
\begin{equation}
\left.  g_{s}\left(  x\right)  \right\vert _{s=0}=1,
\end{equation}
so that the following weighted moment%
\begin{equation}
m_{n}\left(  s\right)  =\int_{-\infty}^{\infty}g_{s}\left(  x\right)  p\left(
x\right)  x^{n}dx \label{mns}%
\end{equation}
converges. The weighted moment recovers the moment when $s=0$:
\begin{equation}
\left.  m_{n}\left(  s\right)  \right\vert _{s=0}=m_{n}.
\end{equation}

If the $n$-th moment does not exist, then $s=0$ must be a singularity of
$m_{n}\left(  s\right)  $. In order to expose the singularity, we expand the
weighted moment $m_{n}\left(  s\right)  $ at $s=0$:
\begin{equation}
m_{n}\left(  s\right)  =%
{\displaystyle\sum_{l=-N}^{\infty}}
a_{l}s^{l}=\frac{a_{-N}}{s^{N}}+\frac{a_{-N+1}}{s^{N+1}}+\cdots+a_{0}%
+a_{1}s+a_{2}s^{2}+\cdots. \label{expansion}%
\end{equation}
The divergence appears in negative-power terms.

In order to obtain a finite $n$-th moment $m_{n}$, we use the minimal
subtraction in quantum field theory. When taking $s=0$ to recover the $n$-th
moment, the negative-power term diverges, and the positive-power term
vanishes. Dropping the divergent negative-power terms, we arrive at a
renormalized moment,
\begin{equation}
m_{n}^{\text{ren}}=a_{0} \label{mnrena0}%
\end{equation}
which is the zeroth-order term in the expansion (\ref{expansion}).

The weighted moment scheme introduces a parameter $s$ so that the moment is a
function of the parameter and the divergence is a singularity of this function.

\subsection{Cut-off scheme \label{IudR}}

Some divergent integrals may become finite after cutting off the integral
limit. For example, when the integral in the moment (\ref{moment}) diverges at
the upper limit, replacing the upper limit with a finite value will give a
cut-off-dependent finite moment. Cutting off the upper limit of an integral is
equivalent to adding a step function as the weight function:%

\begin{equation}
g_{\Lambda}\left(  x\right)  =\theta\left(  \Lambda-x\right)  .
\end{equation}
Then by Eq. (\ref{mns}) we arrive at a cutting-off moment%

\begin{align}
m_{n}\left(  \Lambda\right)   &  =\int_{-\infty}^{\infty}\theta\left(
\Lambda-x\right)  p\left(  x\right)  x^{n}dx\nonumber\\
&  =\int_{-\infty}^{\Lambda}p\left(  x\right)  x^{n}dx.
\end{align}
In this treatment, by introducing a cut-off parameter $\Lambda$, we obtain a
finite cut-off-dependent moment. The same treatment also applies to the cases
where the divergence appears at lower limits and at both upper and lower
limits. The cut-off scheme is essentially a special case of the weighted
moment scheme. Thus, the weighted moment method in the previous section
(\ref{WFR}) is also applicable. If the moment diverges at $\infty$, then
$\Lambda\rightarrow\infty$ is a singular point of $m_{n}\left(  \Lambda
\right)  $. We expand $m_{n}\left(  \Lambda\right)  $ at $\Lambda
\rightarrow\infty$, or equivalently, at $1/\Lambda=0$. The zero-power term of
the expansion of $m_{n}\left(  \Lambda\right)  $, as shown in Eq.
(\ref{mnrena0}), is the renormalized moment. A similar treatment has been used
to deal with the divergence in the scattering of quantum mechanics and proved
to be valid \cite{li2016scattering}.

\subsection{Characteristic function scheme \label{chm}}

The characteristic function is defined as
\begin{equation}
f\left(  k\right)  =\int_{-\infty}^{\infty}p\left(  x\right)  e^{ikx}dx,
\label{chfunction}%
\end{equation}
which is the Fourier transform of the probability density function $p\left(
x\right)  $ \cite{johnson1995continuous}. In this section, we suggest a
renormalization scheme based on the characteristic function. The
characteristic function scheme can also be viewed as a special weight function
scheme with the weight function,
\begin{equation}
g_{k}\left(  x\right)  =e^{ikx}. \label{gsx}%
\end{equation}

The weighted $n$-th moment with the weight function (\ref{gsx}) is
\begin{equation}
m_{n}\left(  k\right)  =\int_{-\infty}^{\infty}p\left(  x\right)  x^{n}%
e^{ikx}dx. \label{mn(s)}%
\end{equation}
It can be seen that the moment (\ref{mn(s)}) is a Fourier transform of the
function $p\left(  x\right)  x^{n}$. This allows us to perform analytical
continuation with the help of the Fourier transform. For example, a
Fourier-type integral $\int_{-\infty}^{\infty}\frac{1}{x}e^{ikx}dx$ is not
integrable, but by analytical continuation the Fourier transform of $\frac
{1}{x}$\ is $i\operatorname*{sgn}\left(  k\right)  \sqrt{\frac{\pi}{2}}$.

In particular, the zeroth moment, corresponding to $n=0$ in Eq. (\ref{mn(s)}%
),
\begin{equation}
m_{0}\left(  k\right)  =\int_{-\infty}^{\infty}p\left(  x\right)
e^{ikx}dx\equiv f\left(  k\right)
\end{equation}
is just the characteristic function.

The characteristic function also plays the role of the moment generating
function, i.e., the $n$-th moment can be expressed as the derivatives of the
characteristic function, or the coefficient of the expansion of the
characteristic function:%
\begin{equation}
m_{n}^{\text{ren}}=\left(  -i\right)  ^{n}\left.  \frac{d^{n}f\left(
k\right)  }{dk^{n}}\right\vert _{k=0}. \label{mndfdk}%
\end{equation}
For the distribution that has no moment, the Fourier transform of\ the
probability density function $p\left(  x\right)  $\ is equivalent to
performing the analytical continuation. Then by Eq. (\ref{mndfdk}), we obtain
a renormalized moment.

If a distribution does not have moments, the derivative of the characteristic
function in Eq. (\ref{mndfdk}) is not well-defined. We need to deal with such
derivatives by analytical continuation. To do this, we use an integral
representation of the derivation \cite{herrmann2018fractional}:
\begin{equation}
\frac{d^{n}f\left(  k\right)  }{dk^{n}}=\frac{1}{\Gamma(-n)}\int_{0}%
^{k}(k-s)^{-n-1}f(s)ds. \label{mnintkdfdkn}%
\end{equation}
Then by Eqs. (\ref{mndfdk}) and (\ref{mnintkdfdkn}), we have
\begin{equation}
m_{n}^{\text{ren}}=\left(  -i\right)  ^{n}\left.  \frac{1}{\Gamma(-n)}\int%
_{0}^{k}(k-s)^{-n-1}f(s)ds\right\vert _{k=0}. \label{mnintk}%
\end{equation}

\subsection{Mellin transform scheme \label{Mellin}}

In this section, we establish a Mellin transform scheme to renormalize power
moments. In this approach, the analytical continuation is performed through
the Mellin transform.

\subsubsection{Characteristic function approach \label{MellinC}}

According to Eq. (\ref{mnintk}), the $n$-th moment is generated by the
characteristic function. We rewrite Eq. (\ref{mnintk}) as
\begin{align}
m_{n}^{\text{ren}}  &  =-\left(  -i\right)  ^{n}\left.  \frac{1}{\Gamma
(-n)}\int_{k}^{\infty}\left(  k-s\right)  ^{-n-1}f\left(  s\right)
ds\right\vert _{k=0}\nonumber\\
&  =-\left(  -i\right)  ^{n}\left(  -1\right)  ^{-n-1}\frac{1}{\Gamma(-n)}%
\int_{0}^{\infty}s^{-n-1}f\left(  s\right)  ds\nonumber\\
&  =\frac{i^{n}}{\Gamma(-n)}\int_{0}^{\infty}s^{-n-1}f\left(  s\right)  ds.
\label{mrenmellinpre}%
\end{align}
By the Mellin transform \cite{brychkov2018handbook},
\begin{equation}
\mathcal{M}_{\sigma}\left[  g\left(  x\right)  \right]  =\int_{0}^{\infty
}g\left(  x\right)  x^{\sigma-1}dx, \label{MLT}%
\end{equation}
we express the $n$-th moment (\ref{mrenmellinpre}) as the Mellin transform of
the characteristic function:
\begin{equation}
m_{n}^{\text{ren}}=\frac{i^{n}}{\Gamma(-n)}\mathcal{M}_{-n}\left[  f\left(
k\right)  \right]  . \label{powermellin}%
\end{equation}

\subsubsection{Density function approach \label{MellinD}}

We can also express the $n$-th moment as the Mellin transform of the
probability density function $p\left(  x\right)  $:
\begin{align}
m_{n}  &  =\int_{-\infty}^{\infty}p\left(  x\right)  x^{n}dx\nonumber\\
&  =\int_{-\infty}^{0}p\left(  x\right)  x^{n}dx+\int_{0}^{\infty}p\left(
x\right)  x^{n}dx\nonumber\\
&  =\int_{\infty}^{0}p\left(  -x\right)  \left(  -x\right)  ^{n}d\left(
-x\right)  +\int_{0}^{\infty}p\left(  x\right)  x^{n}dx\nonumber\\
&  =\left(  -1\right)  ^{n}\int_{0}^{\infty}p\left(  -x\right)  x^{n+1-1}%
dx+\int_{0}^{\infty}p\left(  x\right)  x^{n}dx.
\end{align}
Compared with the definition of the Mellin transform (\ref{MLT}), we have%
\begin{equation}
m_{n}^{\text{ren}}=\left(  -1\right)  ^{n}\mathcal{M}_{n+1}\left[  p\left(
-x\right)  \right]  +\mathcal{M}_{n+1}\left[  p\left(  x\right)  \right]  .
\label{MellinDF}%
\end{equation}
For $x>0$, we have
\begin{equation}
m_{n}^{\text{ren}}=\mathcal{M}_{n+1}\left[  p\left(  x\right)  \right]  .
\end{equation}

\subsection{Power-logarithmic moment scheme \label{LogR}}

In this scheme we introduce the power-logarithmic moment:
\begin{equation}
M_{n,m}=\int_{-\infty}^{\infty}p\left(  x\right)  x^{n}\ln^{m}xdx.
\end{equation}
If the $n$-th power moment%
\begin{equation}
m_{n}=\int_{-\infty}^{\infty}p\left(  x\right)  x^{n}dx \label{mnpxn}%
\end{equation}
is divergent, we can first work out the following power-logarithmic moment
\begin{equation}
M_{n-1,m}=\int_{-\infty}^{\infty}p\left(  x\right)  x^{n-1}\ln^{m}xdx,
\label{Mn-1m}%
\end{equation}
and then achieve the power moment from the power-logarithmic moment by their
relation given below. This approach applies to the distribution\ with which
$p\left(  x\right)  x^{n}$ is not integrable, but $p\left(  x\right)
x^{n-1}\ln^{m}x$ is integrable.

Generally speaking, compared with $p\left(  x\right)  x^{n}$, $p\left(
x\right)  x^{n-1}\ln^{m}x$ is easier to satisfy the integrability condition.
Since $x$ diverges faster than $\ln x$, $\left[  p\left(  x\right)
x^{n-1}\right]  \ln x$ is easier to satisfy the integrability condition than
$\left[  p\left(  x\right)  x^{n-1}\right]  x=p\left(  x\right)  x^{n}$. Thus,
when the $n$-th moment (\ref{mnpxn}) does not exist, the integral
$\int_{-\infty}^{\infty}p\left(  x\right)  x^{n-1}\ln xdx$, even the integral
$\int_{-\infty}^{\infty}p\left(  x\right)  x^{n-1}\ln^{m}xdx$, may possibly
satisfy the integrability condition. In the power-logarithmic moment scheme,
if $\int_{-\infty}^{\infty}p\left(  x\right)  x^{n-1}\ln^{m}xdx$ is
integrable, instead of the integral in (\ref{mnpxn}), we turn to work out\ the
power-logarithmic moment (\ref{Mn-1m}), and then obtain the $n$-th moment
$m_{n}$ from the power-logarithmic moment $M_{n-1,m}$.

A probability density function $p\left(  x\right)  $, with which $p\left(
x\right)  x^{n}$ is not integrable but $p\left(  x\right)  x^{n-1}\ln^{m}x$ is
integrable, has no power moment (\ref{mnpxn}) but has power-logarithmic moment
(\ref{Mn-1m}). By%
\begin{equation}
\ln^{m}x=\frac{1}{x^{n-1}}\frac{d^{m}x^{n-1}}{dn^{m}} \label{Dlog}%
\end{equation}
we have%
\begin{align}
M_{n-1,m}  &  =\frac{d^{m}}{dn^{m}}\int_{-\infty}^{\infty}p\left(  x\right)
x^{n-1}dx\nonumber\\
&  =\frac{d^{m}}{dn^{m}}m_{n-1}. \label{Mn-1}%
\end{align}
Thus, when the power-logarithmic moment $M_{n-1,m}$ is known, the $n$-th
moment, by Eq. (\ref{Mn-1}), can be solved from the following differential
equation:%
\begin{equation}
\frac{d^{m}}{dn^{m}}m_{n}=M_{n,m}.
\end{equation}

In particular, for $m=1$, denoting $M_{n,1}\equiv M_{n}$, we have%
\begin{equation}
\frac{d}{dn}m_{n}=M_{n}. \label{Mn}%
\end{equation}
The $n$-th moment, by solving the differential equation (\ref{Mn}), reads
\begin{equation}
m_{n}^{\text{ren}}=\int M_{n}dn. \label{mnMn}%
\end{equation}
In this way, we can obtain the renormalized $n$-th moment $m_{n}^{\text{ren}}$
if the power-logarithmic moment $M_{n}$ exists.

\section{Renormalized power moment: examples \label{RenEx}}

In this section, we use the renormalization scheme suggested in section
\ref{powerM} to calculate the renormalized $n$-th power moment for the Cauchy
distribution, the Levy distribution, the $q$-exponential distribution, and the
$q$-Gaussian distribution. In order to verify that the renormalization
treatment is scheme-independent, for each distribution, we use more than one
renormalization scheme.

\subsection{Cauchy distribution \label{Cauchy}}

The probability density function of the Cauchy distribution is
\cite{johnson1995continuous}%
\begin{equation}
p\left(  x\right)  =\frac{1}{\pi\left(  1+x^{2}\right)  }. \label{Cauchypdf}%
\end{equation}
The $n$-th moment of the Cauchy distribution, by definition (\ref{norder}), is%
\begin{equation}
m_{n}=\int_{-\infty}^{\infty}\frac{1}{\pi\left(  1+x^{2}\right)  }x^{n}dx.
\label{Cauchymn}%
\end{equation}
The $n$-th moment of the Cauchy distribution does not exist, for the integral
diverges when $n=1,2,3,\cdots$. Next, we use different renormalization
treatments to renormalize the divergent moments.

\subsubsection{Zeta function method}

We first use the zeta function method in section \ref{ZetaR} to calculate the
renormalized $n$-th moment.

The spectral zeta function of the Cauchy distribution, by Eq. (\ref{zetaint}),
is%
\begin{align}
\zeta\left(  s\right)   &  =\int_{-\infty}^{\infty}\frac{x^{-s}}{\pi\left(
1+x^{2}\right)  }dx\nonumber\\
&  =e^{-is\pi/2}. \label{CauchyZeta}%
\end{align}
By Eq. (\ref{mnzeta}), the $n$-th moment of the Cauchy distribution is%
\begin{equation}
m_{n}^{\text{ren}}=\zeta\left(  -n\right)  =e^{in\pi/2}. \label{Cauchymnr}%
\end{equation}
We list the first several moments:
\begin{equation}
m_{1}^{\text{ren}}=i,\text{ }m_{2}^{\text{ren}}=-1,\text{ }m_{3}^{\text{ren}%
}=-i,\text{ }m_{4}^{\text{ren}}=1,\text{ }\cdots.
\end{equation}
It can be seen that for integer-order moments, we have%
\begin{align}
m_{2n}^{\text{ren}}  &  =\left(  -1\right)  ^{n},\text{ }\ \text{\ }%
n=\text{even,}\nonumber\\
m_{2n-1}^{\text{ren}}  &  =\left(  -1\right)  ^{n-1}i,\text{ \ }n=\text{odd.}%
\end{align}

It can be seen that some of the renormalized moments are less than zero, and
some are pure imaginary numbers.

\subsubsection{Subtraction method}

We use the subtraction method in section \ref{remove} to calculate the
renormalized $n$-th moment.

The $n$-th moment of the Cauchy distribution (\ref{Cauchymn}) diverges at
$x\rightarrow\pm\infty$. Expand the integrand in Eq. (\ref{Cauchymn}) at
$x\rightarrow-\infty$ and $x\rightarrow\infty$:
\begin{equation}
\frac{x^{n}}{\pi\left(  1+x^{2}\right)  }=-\frac{1}{\pi}\sum_{j=1}^{\infty
}\cos\left(  \frac{j\pi}{2}\right)  \frac{1}{x^{j-n}}.
\end{equation}
The divergence is caused by the term $j-n\leq1$, i.e.,
\begin{equation}
j\leq n+1.
\end{equation}
By Eq. (\ref{mreninf+-}), we have
\begin{align}
m_{n}^{\text{ren}}  &  =\int_{-\infty}^{-1}\left[  \frac{x^{n}}{\pi\left(
1+x^{2}\right)  }+\frac{1}{\pi}\sum_{j=1}^{j\leq n+1}\cos\left(  \frac{j\pi
}{2}\right)  \frac{1}{x^{j-n}}\right]  dx+i\pi c_{-1}-\sum_{j=1}^{j\leq
n+1}\cos\left(  \frac{j\pi}{2}\right)  \frac{\left(  -1\right)  ^{j-n-1}%
}{1-\left(  j-n\right)  }\nonumber\\
&  +\int_{-1}^{1}\frac{x^{n}}{\pi\left(  1+x^{2}\right)  }dx+\int_{1}^{\infty
}\left[  \frac{x^{n}}{\pi\left(  1+x^{2}\right)  }+\frac{1}{\pi}\sum
_{j=1}^{j\leq n+1}\cos\left(  \frac{j\pi}{2}\right)  \frac{1}{x^{j-n}}\right]
dx\nonumber\\
&  +\sum_{j=1}^{j\leq n+1}\cos\left(  \frac{j\pi}{2}\right)  \frac
{1}{1-\left(  j-n\right)  }.
\end{align}
For example, when $n=1$, the first-order moment is%
\begin{align}
m_{1}^{\text{ren}}  &  =\int_{-\infty}^{-1}\left[  \frac{x}{\pi\left(
1+x^{2}\right)  }-\frac{1}{\pi x}\right]  dx+i\pi c_{-1}+\int_{-1}^{1}\frac
{x}{\pi\left(  1+x^{2}\right)  }dx+\int_{1}^{\infty}\left[  \frac{x^{n}}%
{\pi\left(  1+x^{2}\right)  }-\frac{1}{\pi x}\right]  dx\nonumber\\
&  =i,
\end{align}
where the coefficient of $\frac{1}{x}$ in the expansion is $c_{-1}=\frac
{1}{\pi}m$.

\subsubsection{Weight function method (1)}

We use the weighted moment method in section \ref{WFR} to calculate the
renormalized $n$-th moment. In order to exemplify that the renormalized moment
is independent of the choice of weight functions, we choose two weight functions.

Taking the weight function as
\begin{equation}
g_{s}\left(  x\right)  =e^{-sx}%
\end{equation}
gives the weighted moment
\begin{align}
m_{n}\left(  s\right)   &  =\int_{-\infty}^{\infty}\frac{1}{\pi\left(
1+x^{2}\right)  }x^{n}e^{-sx}dx\nonumber\\
&  =e^{i\pi n/2}\left[  \cos s+\operatorname*{sgn}(s)\sin s-\left(
\text{$\operatorname*{sgn}$}(s)+i\right)  \frac{i^{-n}s^{-n+1}}{\Gamma
(2-n)}\,_{1}F_{2}\left(  1;1-\frac{n}{2},\frac{3}{2}-\frac{n}{2};-\frac{s^{2}%
}{4}\right)  \right]  .
\end{align}
\bigskip According to Eq. (\ref{mnrena0}), we expand the weighted moment
$m_{n}\left(  s\right)  $ at $s=0$:
\begin{equation}
e^{in\pi/2}\left[  1+s-\frac{(1+i)s(is)^{-n}}{\Gamma(2-n)}+\cdots\right]  ,
\end{equation}
where the zero-power term is the renormalized moment
\begin{equation}
m_{n}^{\text{ren}}\left(  0\right)  =e^{in\pi/2}.
\end{equation}

\subsubsection{Weight function method (2)}

Taking another the weight function,%
\begin{equation}
g_{s}\left(  x\right)  =x^{s-1},
\end{equation}
the weighted moment is%
\begin{align}
m_{n}\left(  s\right)   &  =\int_{-\infty}^{\infty}\frac{1}{\pi\left(
1+x^{2}\right)  }x^{n}x^{s-1}dx\nonumber\\
&  =\frac{1}{2}\left[  (-1)^{n+s+1}+1\right]  \csc\left(  \frac{1}{2}%
\pi(n+s)\right)  .
\end{align}
According to Eq. (\ref{mnrena0}), we expand the weighted moment $m_{n}\left(
s\right)  $ at $s=0$:%
\begin{equation}
e^{in\pi/2}\left[  1+\frac{i\pi(s-1)}{2}+\cdots\right]  ,
\end{equation}
where the zero-power term is the renormalized moment%
\begin{equation}
m_{n}^{\text{ren}}=e^{in\pi/2}.
\end{equation}

\subsubsection{Cut-off method}

We use the cut-off method in section \ref{IudR} to calculate the renormalized
$n$-th moment.

Cutting off the upper and lower integration limit of the $n$-th moment
(\ref{Cauchymn}) gives
\begin{equation}
m_{n}\left(  \Lambda_{1},\Lambda_{2}\right)  =\int_{-\Lambda_{2}}^{\Lambda
_{1}}\frac{1}{\pi\left(  1+x^{2}\right)  }x^{n}dx. \label{cuachycut}%
\end{equation}
When $\Lambda_{1}\rightarrow\infty$ and $\Lambda_{2}\rightarrow\infty$, we
obtain the $n$-th moment: $m_{n}=m_{n}\left(  \infty,\infty\right)  $. Working
out the integral in Eq. (\ref{cuachycut}) gives
\begin{equation}
m_{n}\left(  \Lambda_{1},\Lambda_{2}\right)  =\frac{\Lambda_{1}^{n+1}%
\,_{2}F_{1}\left(  1,\frac{n+1}{2};\frac{n+3}{2};-\Lambda_{1}^{2}\right)
-\left(  -\Lambda_{2}\right)  ^{n+1}\,_{2}F_{1}\left(  1,\frac{n+1}{2}%
;\frac{n+3}{2};-\Lambda_{2}^{2}\right)  }{\pi(n+1)}.
\end{equation}
In order to obtain the renormalized $n$-th moment, we expand $m_{n}\left(
\Lambda\right)  $ at $\Lambda_{1},$ $\Lambda_{2}\rightarrow\infty$:
\begin{equation}
m_{n}\left(  \Lambda\right)  =e^{in\pi/2}+\frac{\Lambda_{1}^{n-1}}{\pi
(n-1)}-\frac{\left(  -\Lambda_{2}\right)  ^{n-1}}{\pi(n-1)}+\cdots.
\end{equation}
Dropping the divergent terms, we obtain the renormalized $n$-th moment,
\begin{equation}
m_{n}^{\text{ren}}=e^{in\pi/2}.
\end{equation}

\subsubsection{Characteristic function method}

We use the characteristic function method in section \ref{chm} to calculate
the renormalized $n$-th moment.

The characteristic function of the Cauchy distribution by Eq.
(\ref{chfunction}) is%
\begin{equation}
f\left(  k\right)  =\int_{-\infty}^{\infty}p\left(  x\right)  e^{ikx}%
dx=e^{-\left\vert k\right\vert }. \label{CauchyC}%
\end{equation}

By the integral representation of the derivation, Eq. (\ref{mnintkdfdkn}), we
obtain
\begin{equation}
\frac{d^{n}f\left(  k\right)  }{dk^{n}}=e^{-k+i\pi n}\frac{\Gamma
(-n)-\Gamma(-n,-k)}{\Gamma(-n)}.
\end{equation}
By Eq. (\ref{mndfdk}), the renormalized $n$-th moment is
\begin{equation}
m_{n}^{\text{ren}}=e^{\frac{i\pi n}{2}}\frac{\Gamma(-n)-\Gamma(-n,0)}%
{\Gamma(-n)}.
\end{equation}
When $n$ is a positive integer, the renormalized $n$-th moment is%
\[
m_{n}^{\text{ren}}=e^{in\pi/2}.
\]

\subsubsection{Mellin transform scheme: characteristic function method}

We use the Mellin transform method for the characteristic function in section
\ref{MellinC} to calculate the renormalized $n$-th moment.

Substituting the characteristic function of the Cauchy distribution
(\ref{CauchyC}) into\ Eq. (\ref{powermellin}), we obtain%
\begin{align}
m_{n}^{\text{ren}}  &  =\frac{i^{n}}{\Gamma(-n)}\mathcal{M}_{-n}\left[
e^{-\left\vert k\right\vert }\right]  =\frac{i^{n}}{\Gamma(-n)}\times
\Gamma(-n)=i^{n}\nonumber\\
&  =e^{in\pi/2}.
\end{align}

\subsubsection{Mellin transform scheme: probability density function}

We use the Mellin transform method for the probability density function in
section \ref{MellinD} to calculate the renormalized $n$-th moment.

Substituting the probability density function of the Cauchy distribution
(\ref{Cauchypdf}) into Eq. (\ref{MellinDF}), we obtain%
\begin{align}
m_{n}^{\text{ren}}  &  =\left(  -1\right)  ^{n}\mathcal{M}_{n+1}\left[
p\left(  -x\right)  \right]  +\mathcal{M}_{n+1}\left[  p\left(  x\right)
\right] \nonumber\\
&  =\left(  -1\right)  ^{n}\left(  \frac{1}{2}\sec\frac{n\pi}{2}\right)
+\frac{1}{2}\sec\frac{n\pi}{2}\nonumber\\
&  =e^{in\pi/2}.
\end{align}

\subsubsection{Power-logarithmic moment method}

The power-logarithmic moment of the Cauchy distribution, by Eq. (\ref{Mn-1m}),
is
\begin{align}
M_{n-1}  &  =\int_{-\infty}^{\infty}\frac{1}{\pi\left(  1+x^{2}\right)
}x^{n-1}\ln xdx\nonumber\\
&  =\frac{\pi}{2}e^{in\pi/2}.
\end{align}
By the relation between the $n$-th\ moment $m_{n}$ and the power-logarithmic
moment $M_{n-1}$, Eq. (\ref{mnMn}), we have
\begin{equation}
m_{n}^{\text{ren}}=e^{in\pi/2}.
\end{equation}
The renormalized $n$-th moment given by various renormalization schemes is the
same. This verifies that the renormalization treatment is scheme-independent.

\subsection{Levy distribution}

The probability density function of the Levy distribution is
\cite{johnson1995continuous}%
\begin{equation}
p\left(  x\right)  =\frac{1}{\sqrt{2\pi}}\frac{1}{x^{3/2}}e^{-1/\left(
2x\right)  },\text{ \ \ }0\leq x<\infty. \label{Levypdf}%
\end{equation}
The $n$-th moment of the Levy distribution, by definition (\ref{norder}), is%
\begin{equation}
m_{n}=\int_{0}^{\infty}\frac{1}{\sqrt{2\pi}}\frac{1}{x^{3/2}}e^{-1/\left(
2x\right)  }x^{n}dx. \label{Levynm}%
\end{equation}

\bigskip

The $n$-th moment of the Levy distribution does not exist, for the integral
diverges when $n=1,2,3,\cdots$. In order to verify that the renormalization
treatment is scheme-independent, we use different renormalization treatments
to renormalize the divergent moments.

\subsubsection{Zeta function method}

We first use the zeta function method in section \ref{ZetaR} to calculate the
renormalized $n$-th moment.

The spectral zeta function of the Levy distribution, by Eq. (\ref{zetaint}),
is%
\begin{align}
\zeta\left(  s\right)   &  =\int_{0}^{\infty}\frac{e^{-1/2x}}{\sqrt{2\pi}%
}\frac{1}{x^{3/2}}x^{-s}dx\nonumber\\
&  =\frac{2^{s}}{\sqrt{\pi}}\Gamma\left(  \frac{1}{2}+s\right)  .
\end{align}

By Eq. (\ref{mnzeta}), the $n$-th moment of the Levy distribution is%
\begin{equation}
m_{n}^{\text{ren}}=\zeta\left(  -n\right)  =\frac{1}{2^{n}\sqrt{\pi}}%
\Gamma\left(  \frac{1}{2}-n\right)  . \label{Levymnr}%
\end{equation}
The first several renormalized $n$-th moments can be obtained directly,
\begin{equation}
m_{1}^{\text{ren}}=-1,\text{ }m_{2}^{\text{ren}}=\frac{1}{3},\text{ }%
m_{3}^{\text{ren}}=-\frac{1}{15},\text{ }m_{4}^{\text{ren}}=\frac{1}%
{105},\cdots.
\end{equation}

\subsubsection{Subtraction method}

We use the subtraction method in section \ref{removeinf} to calculate the
renormalized $n$-th moment of the Levy distribution.

The $n$-th moment of the Levy distribution (\ref{Levynm}) diverges at
$x\rightarrow\infty$. Expand the integrand at $x\rightarrow\infty$:%
\begin{equation}
\frac{e^{-\frac{1}{2x}}}{\sqrt{2\pi}}\frac{1}{x^{3/2}}x^{n}=\sum_{j=0}%
^{\infty}\frac{\left(  -1\right)  ^{j}2^{-\frac{1}{2}-j}}{\sqrt{\pi}j!}%
\frac{1}{x^{-n+3/2+j}}.
\end{equation}
The divergence is caused by the term $-n+3/2+j\leq1$, i.e., $j\leq n-\frac
{1}{2}$. Since $j$ is an integer, the divergent term is
\begin{equation}
j<n.
\end{equation}
By Eq. (\ref{mreninf}), we obtain the renormalized $n$-th moment
\begin{align}
m_{n}^{\text{ren}}  &  =\int_{0}^{1}\left(  \frac{e^{-\frac{1}{2x}}}%
{\sqrt{2\pi}}\frac{1}{x^{3/2}}x^{n}\right)  dx+\int_{1}^{\infty}\left(
\frac{e^{-\frac{1}{2x}}}{\sqrt{2\pi}}\frac{1}{x^{3/2}}x^{n}-\sum_{j=0}%
^{j<n}\frac{\left(  -1\right)  ^{j}2^{-\frac{1}{2}-j}}{\sqrt{\pi}j!}\frac
{1}{x^{-n+3/2+j}}\right)  dx\nonumber\\
&  -\sum_{j=0}^{j<n}\frac{\left(  -1\right)  ^{j}2^{-\frac{1}{2}-j}}{\sqrt
{\pi}j!}\frac{1}{1-\left(  -n+\frac{3}{2}+j\right)  }.
\end{align}
In particular, the first-order moment is%
\begin{align}
m_{1}^{\text{ren}}  &  =\int_{0}^{1}\left(  \frac{e^{-\frac{1}{2x}}}%
{\sqrt{2\pi}}\frac{1}{x^{3/2}}x\right)  dx+\int_{1}^{\infty}\left(
\frac{e^{-\frac{1}{2x}}}{\sqrt{2\pi}}\frac{1}{x^{3/2}}x-\frac{1}{\sqrt{2\pi}%
}\frac{1}{x^{1/2}}\right)  dx-\frac{\frac{1}{\sqrt{2\pi}}}{1-\frac{1}{2}%
}\nonumber\\
&  =\sqrt{\frac{2}{e\pi}}-\operatorname*{erfc}\left(  \frac{1}{\sqrt{2}%
}\right)  +\left(  \sqrt{\frac{2}{\pi}}-\sqrt{\frac{2}{e\pi}}%
-\operatorname{erf}\left(  \frac{1}{\sqrt{2}}\right)  \right)  -\sqrt{\frac
{2}{\pi}}\nonumber\\
&  =-1.
\end{align}

\subsubsection{Weight function method (1)}

We use the weighted moment method in section \ref{WFR} to calculate the
renormalized $n$-th moment of the Levy distribution. In order to exemplify
that the renormalized moment is independent of the choice of weight functions,
we use two weight functions.

Taking the weight function as%
\begin{equation}
g_{s}\left(  x\right)  =e^{-sx}%
\end{equation}
gives the weighted moment
\begin{align}
m_{n}\left(  s\right)   &  =\int_{0}^{\infty}\frac{e^{-\frac{1}{2x}}}%
{\sqrt{2\pi}}\frac{1}{x^{3/2}}x^{n}e^{-sx}dx\nonumber\\
&  =\frac{1}{\sqrt{\pi}}2^{\frac{3}{4}-\frac{n}{2}}s^{\frac{1}{4}-\frac{n}{2}%
}K_{n-\frac{1}{2}}\left(  \sqrt{2s}\right)  ,
\end{align}
where $K_{n}\left(  z\right)  $ is the modified Bessel function of the second
kind \cite{olver2010nist}. According to Eq. (\ref{mnrena0}), we expand the
weighted moment $m_{n}\left(  s\right)  $ at $s=0$:
\begin{equation}
\frac{1}{2^{n}\sqrt{\pi}}\Gamma\left(  \frac{1}{2}-n\right)  +\frac
{2^{-n}\Gamma\left(  \frac{1}{2}-n\right)  }{\sqrt{\pi}(2n+1)}s+\frac
{\Gamma\left(  n-\frac{1}{2}\right)  }{\sqrt{2\pi}}s^{\frac{1}{2}-n}+\cdots,
\end{equation}
where the zero-order term is the renormalized moment:%
\begin{equation}
m_{n}^{\text{ren}}\left(  0\right)  =\frac{1}{2^{n}\sqrt{\pi}}\Gamma\left(
\frac{1}{2}-n\right)  .
\end{equation}
In particular, the first moment is $m_{1}^{\text{ren}}=-1$.

\subsubsection{Weight function method (2)}

Taking another the weight function,%
\begin{equation}
g_{s}\left(  x\right)  =x^{s-1},
\end{equation}
the weight moment is%
\begin{align}
m_{n}\left(  s\right)   &  =\int_{0}^{\infty}\frac{e^{-\frac{1}{2x}}}%
{\sqrt{2\pi}}\frac{1}{x^{3/2}}x^{n}x^{s-1}dx\nonumber\\
&  =\frac{2^{-n-s+1}}{\sqrt{\pi}}\Gamma\left(  -n-s+\frac{3}{2}\right)  .
\end{align}
According to (\ref{mnrena0}), we expand the weighted moment $m_{n}\left(
s\right)  $ at $s=0$:%
\begin{equation}
\frac{1}{2^{n}\sqrt{\pi}}\Gamma\left(  \frac{1}{2}-n\right)  -\frac{1}%
{2^{n}\sqrt{\pi}}(s-1)\Gamma\left(  \frac{1}{2}-n\right)  \left(  \psi
^{(0)}\left(  \frac{1}{2}-n\right)  +\ln2\right)  +\cdots,
\end{equation}
where the zero-order term is the renormalization $n$-th moment:%
\begin{equation}
m_{n}^{\text{ren}}\left(  1\right)  =\frac{1}{2^{n}\sqrt{\pi}}\Gamma\left(
\frac{1}{2}-n\right)  .
\end{equation}

\subsubsection{Cut-off method}

We use the cut-off method in section \ref{IudR} to calculate the renormalized
$n$-th moment.

Cutting off the upper integration limit of the $n$-th moment gives%
\begin{equation}
m_{n}\left(  \Lambda\right)  =\int_{0}^{\Lambda}\frac{e^{-1/\left(  2x\right)
}}{\sqrt{2\pi}}\frac{1}{x^{3/2}}x^{n}dx. \label{levycut}%
\end{equation}
When $\Lambda\rightarrow\infty$, we obtain the $n$-th moment: $m_{n}%
=m_{n}\left(  \infty\right)  $. Working out the integral in Eq. (\ref{levycut}%
) gives%
\begin{equation}
m_{n}\left(  \Lambda\right)  =\frac{1}{2^{n}\sqrt{\pi}}\Gamma\left(  \frac
{1}{2}-n,\frac{1}{2\Lambda}\right)  .
\end{equation}
In order to obtain the renormalized $n$-th moment, we expand $m_{n}\left(
\Lambda\right)  $ at $\Lambda\rightarrow\infty$:%
\begin{equation}
m_{n}\left(  \Lambda\right)  =\frac{1}{2^{n}\sqrt{\pi}}\Gamma\left(  \frac
{1}{2}-n\right)  +\frac{\sqrt{2}}{\sqrt{\pi}\left(  2n-1\right)  }%
\Lambda^{n-1/2}+\cdots.
\end{equation}
The zero-order term is the renormalized $n$-th moment,%
\begin{equation}
m_{n}^{\text{ren}}=\frac{1}{2^{n}\sqrt{\pi}}\Gamma\left(  \frac{1}%
{2}-n\right)  .
\end{equation}

\subsubsection{Characteristic function method}

We use the characteristic function method in section \ref{chm} to calculate
the renormalized $n$-th moment of the Levy distribution.

The characteristic function of the Levy distribution by Eq. (\ref{chfunction})
is%
\begin{equation}
f\left(  k\right)  =e^{-\sqrt{\left\vert k\right\vert }}\left(  \cos
\sqrt{\left\vert k\right\vert }+i\operatorname*{sgn}(k)\sin\sqrt{\left\vert
k\right\vert }\right)  . \label{LevyC}%
\end{equation}
By the integral representation of the derivation, Eq. (\ref{mnintkdfdkn}), we
obtain%
\begin{align}
\frac{d^{n}f\left(  k\right)  }{dk^{n}}  &  =\pi^{3/2}2^{n-\frac{5}{2}}%
k^{-n}\left[  8_{1}\tilde{F}_{4}\left(  1;\frac{1}{4},\frac{3}{4},\frac{1}%
{2}-\frac{n}{2},1-\frac{n}{2};-\frac{k^{2}}{64}\right)  \right. \nonumber\\
&  -\left.  ik\text{ }_{1}\tilde{F}_{4}\left(  1;\frac{3}{4},\frac{5}%
{4},1-\frac{n}{2},\frac{3}{2}-\frac{n}{2};-\frac{k^{2}}{64}\right)  \right]
\nonumber\\
&  +(-1)^{5/8}\sqrt{\pi}2^{-\frac{n}{2}-\frac{1}{4}}e^{-\frac{3}{4}i\pi
n}k^{\frac{1}{4}-\frac{n}{2}}\left[  \operatorname*{bei}\nolimits_{\frac{1}%
{2}-n}\left(  \sqrt{2k}\right)  +i\operatorname*{ber}\nolimits_{\frac{1}{2}%
-n}\left(  \sqrt{2k}\right)  \right]  , \label{Lcren}%
\end{align}
where $_{p}\tilde{F}_{q}\left(  a;b;z\right)  $ is regularized hypergeometric
function, given by $_{p}F_{q}\left(  a;b;z\right)  /\left(  \Gamma
(b_{1})\cdots\Gamma(b_{q})\right)  $, $\operatorname*{bei}_{\nu}\left(
z\right)  $ and $\operatorname*{ber}\nolimits_{\nu}\left(  z\right)  $ are the
Kelvin function and $\operatorname*{bei}_{\nu}\left(  z\right)
+i\operatorname*{ber}\nolimits_{\nu}\left(  z\right)  =J_{\nu}\left(
xe^{\frac{3\pi i}{4}}\right)  $, and $J_{\nu}\left(  x\right)  $ is a Bessel
function of the first kind \cite{olver2010nist}.

By Eqs. (\ref{mndfdk}) and (\ref{Lcren}), we have%
\begin{align}
m_{1}^{\text{ren}}  &  =\left(  -i\right)  \left.  \frac{df\left(  k\right)
}{dk}\right\vert _{k=0}=-1,\nonumber\\
m_{2}^{\text{ren}}  &  =\left(  -i\right)  ^{2}\left.  \frac{d^{2}f\left(
k\right)  }{dk^{2}}\right\vert _{k=0}=\frac{1}{3},\nonumber\\
&  \cdots.
\end{align}

\subsubsection{Mellin transform method: characteristic function}

We use the Mellin transform method of the characteristic function in section
\ref{MellinC} to calculate the renormalized $n$-th moment of the Levy distribution.

Substituting the characteristic function of the Levy distribution, Eq.
(\ref{LevyC}), into Eq. (\ref{powermellin}), we obtain%
\begin{align}
m_{n}^{\text{ren}}  &  =\frac{i^{n}}{\Gamma(-n)}\mathcal{M}_{-n}\left[
f\left(  k\right)  \right] \nonumber\\
&  =\frac{i^{n}}{\Gamma(-n)}2^{n+1}\left(  -i\right)  ^{n}\Gamma\left(
-2n\right) \nonumber\\
&  =\frac{1}{2^{n}\sqrt{\pi}}\Gamma\left(  \frac{1}{2}-n\right)  .
\end{align}

\subsubsection{Mellin transform method: probability density function}

We use the Mellin transform method of the probability density function in
section \ref{MellinD} to calculate the renormalized $n$-th moment of the Levy distribution.

Substituting the probability density function of the Levy distribution, Eq.
(\ref{Levypdf}) into Eq. (\ref{MellinDF}), we obtain%
\begin{align}
m_{n}^{\text{ren}}  &  =\mathcal{M}_{n+1}\left[  p\left(  x\right)  \right]
\nonumber\\
&  =\frac{1}{2^{n}\sqrt{\pi}}\Gamma\left(  \frac{1}{2}-n\right)  .
\end{align}

\subsubsection{Power-logarithmic moment method}

We use the power-logarithmic moment method in section \ref{LogR} to calculate
the renormalized $n$-th moment of the Levy distribution. The power-logarithmic
moment of the Levy distribution, by Eq. (\ref{Mn-1m}), is%
\begin{align}
M_{n-1}  &  =\int_{-\infty}^{\infty}\frac{e^{-\frac{1}{2x}}}{\sqrt{2\pi}}%
\frac{1}{x^{3/2}}x^{n-1}\ln xdx\nonumber\\
&  =-\frac{1}{\sqrt{\pi}}2^{1-n}\Gamma\left(  \frac{3}{2}-n\right)  \left(
\psi^{(0)}\left(  \frac{3}{2}-n\right)  +\ln2\right)  , \label{Mn-1Levy}%
\end{align}
where $\psi^{(n)}\left(  z\right)  $ is the $n$-th derivative of the digamma
function, $\psi^{(n)}\left(  z\right)  =\frac{d^{n}\psi\left(  z\right)
}{dz^{n}}$, and $\psi\left(  z\right)  $ is the logarithmic derivative of the
gamma function, $\psi\left(  z\right)  =\frac{\Gamma^{\prime}\left(  z\right)
}{\Gamma\left(  z\right)  }$ \cite{olver2010nist}. Substituting Eq.
(\ref{Mn-1Levy}) into the relation between the $n$-th $m_{n}$\ moment and the
power-logarithmic moment $M_{n-1}$, Eq. (\ref{mnMn}), we have
\begin{equation}
m_{n}^{\text{ren}}=\frac{1}{2^{n}\sqrt{\pi}}\Gamma\left(  \frac{1}%
{2}-n\right)  .
\end{equation}

The renormalized $n$-th moment given by the different renormalization schemes
is the same. This verifies that the renormalization treatment is scheme-independent.

\subsection{$q$-exponential distribution}

The $q$-exponential distribution has power moments for some values of the
parameter $q$ and no power moments for others. We only need to renormalize the
case having no power moment. The self-consistency requires that the
renormalized power moment be valid for the $q$ with which the power moment is
well-defined. That is, the renormalized $n$-th moment is valid for all $q$.

The probability density function of the $q$-exponential distribution is
\cite{kampfer2021impact}%
\begin{equation}
p\left(  x\right)  =\left\{
\begin{array}
[c]{c}%
\lambda(2-q)(1-\lambda(1-q)x)^{\frac{1}{1-q}},\text{ \ \ }q<1,\text{ }0\leq
x<\frac{1}{\lambda(1-q)},\\
\lambda e^{-\lambda x},\text{ \ }q=1,\text{ \ }x\geq0,\\
\lambda(2-q)(1-\lambda(1-q)x)^{\frac{1}{1-q}},\text{ \ }q>1,\text{ \ }x\geq0.
\end{array}
\right.  . \label{q-exp}%
\end{equation}
For $q=1$, the $q$-exponential distribution returns to the exponential distribution.

The $n$-th moment of the $q$-exponential distribution, by definition
(\ref{norder}), is%
\begin{equation}
m_{n}=\int_{0}^{\infty}\lambda(2-q)\left[  1-\lambda(1-q)x\right]  ^{\frac
{1}{1-q}}x^{n}dx. \label{qemn}%
\end{equation}
For $q\leq1$, the $n$-th moment $m_{n}$ exists. For $q>1$, the $n$-th moment
$m_{n}$ exists when $q<\frac{n+2}{n+1}$. For example, only for $q<\frac{3}{2}%
$, there exists the first moment, $m_{1}=\frac{1}{3\lambda-2\lambda q}$; only
for $q<\frac{4}{3}$, there exists the second moment, $m_{2}=\frac{2}%
{\lambda^{2}\left(  6q^{2}-17q+12\right)  }$; only for $q<\frac{5}{4}$, there
exists the third moment, $m_{3}=-\frac{6}{\lambda^{3}\left(  24q^{3}%
-98q^{2}+133q-60\right)  }$; only for $q<\frac{6}{5}$, there exists the fourth
moment, $m_{4}=\frac{24\Gamma\left(  \frac{1}{q-1}-5\right)  }{\lambda
^{4}(q-1)^{4}\Gamma\left(  \frac{1}{q-1}-1\right)  }$, and so on; otherwise
the $n$-th moment does not exist.

Next, we calculate the renormalized $n$-th moment for the case having no power moments.

\subsubsection{Zeta function method}

We first use the zeta function method in section \ref{ZetaR} to calculate the
renormalized $n$-th moment.

The spectral zeta function of the $q$-exponential distribution, by
(\ref{zetaint}), is%
\begin{align}
\zeta\left(  s\right)   &  =\int_{0}^{\infty}\lambda(2-q)(1-\lambda
(1-q)x)^{\frac{1}{1-q}}x^{-s}dx\nonumber\\
&  =\frac{\Gamma(1-s)\Gamma\left(  \frac{-s-q(1-s)+2}{q-1}\right)  }%
{\Gamma\left(  \frac{1}{q-1}-1\right)  \left[  \lambda\left(  q-1\right)
\right]  ^{s}}.
\end{align}
The $n$-th moment of the $q$-exponential distribution, by Eq. (\ref{mnzeta}),
is%
\begin{align}
m_{n}^{\text{ren}}  &  =\zeta\left(  -n\right) \nonumber\\
&  =\frac{\Gamma(n+1)\Gamma\left(  -n+\frac{1}{q-1}-1\right)  }{\Gamma\left(
\frac{1}{q-1}-1\right)  \left[  \lambda(q-1)\right]  ^{n}}. \label{qemnr}%
\end{align}
\qquad

For integer $n$, $m_{n}^{\text{ren}}$ can be written as
\begin{equation}
m_{n}^{\text{ren}}=\left(  -1\right)  ^{n}\frac{n!}{\lambda^{n}}\frac{1}%
{\prod\limits_{l=1}^{n}\left[  \left(  l+1\right)  q-\left(  l+2\right)
\right]  }.
\end{equation}
Thus, the first several moments are
\begin{align}
m_{1}^{\text{ren}}  &  =\zeta\left(  -1\right)  =-\frac{1}{\lambda\left(
2q-3\right)  },\nonumber\\
m_{2}^{\text{ren}}  &  =\zeta\left(  -2\right)  =\frac{2}{\lambda
^{2}(2q-3)(3q-4)},\nonumber\\
m_{3}^{\text{ren}}  &  =\zeta\left(  -3\right)  =-\frac{6}{\lambda
^{3}(2q-3)(3q-4)(4q-5)},\nonumber\\
m_{4}^{\text{ren}}  &  =\zeta\left(  -4\right)  =\frac{24}{\lambda
^{4}(2q-3)(3q-4)(4q-5)(5q-6)}.
\end{align}

\subsubsection{Subtraction method}

We use the subtraction method in section \ref{remove} to calculate the
renormalized $n$-th moment of the $q$-exponential distribution.

The $n$-th moment of the $q$-exponential distribution (\ref{qemn}) diverges at
$x\rightarrow\infty$. Expand the integrand at $x\rightarrow\infty$:
\begin{align}
\lambda(2-q)\left[  1+\lambda(q-1)x\right]  ^{\frac{1}{1-q}}x^{n}  &
=\lambda(2-q)x^{-\frac{1}{q-1}}\left[  \frac{1}{x}+\lambda(q-1)\right]
^{\frac{1}{1-q}}x^{n}\nonumber\\
&  =\lambda\left(  2-q\right)  \sum_{j=0}^{\infty}\left(
\begin{array}
[c]{c}%
\frac{1}{1-q}\\
j
\end{array}
\right)  \left[  \lambda(q-1)\right]  ^{\frac{1}{1-q}-j}\frac{1}%
{x^{j-n+\frac{1}{q-1}}},
\end{align}
where $\left(
\begin{array}
[c]{c}%
m\\
n
\end{array}
\right)  =\frac{\Gamma\left(  m+1\right)  }{\Gamma\left(  n+1\right)
\Gamma\left(  m-n+1\right)  }$\ is the binomial expansion coefficient.
Following section \ref{removeinf}, divergence is caused by the term
$j-n+\frac{1}{q-1}\leq1$, i.e.,
\begin{equation}
j\leq1+n-\frac{1}{q-1}. \label{jint}%
\end{equation}
By Eq. (\ref{mreninf}), we arrive at the renormalized $n$-th moment,
\begin{align}
m_{n}^{\text{ren}}  &  =\int_{0}^{1}\lambda(2-q)\left[  1+\lambda
(q-1)x\right]  ^{\frac{1}{1-q}}x^{n}dx+\int_{1}^{\infty}\left[  \lambda
(2-q)\left[  1+\lambda(q-1)x\right]  ^{\frac{1}{1-q}}x^{n}\right. \nonumber\\
&  -\left.  \lambda\left(  2-q\right)  \sum_{j=0}^{j\leq1+n-\frac{1}{q-1}%
}\left(
\begin{array}
[c]{c}%
\frac{1}{1-q}\\
j
\end{array}
\right)  \left[  \lambda(q-1)\right]  ^{\frac{1}{1-q}-j}\frac{1}%
{x^{j-n+\frac{1}{q-1}}}\right]  dx\nonumber\\
&  -\lambda\left(  2-q\right)  \sum_{j=0}^{j\leq1+n-\frac{1}{q-1}}\left(
\begin{array}
[c]{c}%
\frac{1}{1-q}\\
j
\end{array}
\right)  \frac{\left[  \lambda(q-1)\right]  ^{\frac{1}{1-q}-j}}{1-\left(
j-n+\frac{1}{q-1}\right)  }.
\end{align}
This gives%
\begin{align}
m_{n}^{\text{ren}}  &  =\left(  2-q\right)  \lambda\Gamma\left(  n+1\right)
\text{ }_{2}\widetilde{F}_{1}\left(  n+1,\frac{1}{q-1},n+2,-\lambda\left(
q-1\right)  \right) \nonumber\\
&  +\frac{\left(  -1\right)  ^{-\frac{1}{q-1}}\left(  q-2\right)  }%
{\lambda^{n}\left(  1-q\right)  ^{n}}B\left(  \frac{1}{\lambda\left(
1-q\right)  },\frac{1}{q-1}-n-1,\frac{q-2}{q-1}\right)  , \label{qexp1}%
\end{align}
where $_{2}\widetilde{F}_{1}\left(  \alpha,\beta,\gamma,z\right)  $\ is the
regularized hypergeometric function and $B\left(  \alpha,\beta,z\right)  $ is
the Beta function.

For example, the first-order moment, by (\ref{qexp1}), is%
\begin{equation}
m_{1}^{\text{ren}}=\frac{1}{\left(  3-2q\right)  \lambda}.
\end{equation}

\subsubsection{Weight function method (1)}

We use the weighted moment method in section \ref{WFR} to calculate the
renormalized $n$-th moment of the $q$-exponential distribution. In order to
exemplify that the renormalized moment is independent of the choice of weight
functions, we choose two weight functions.

Taking the weight function as%
\begin{equation}
g_{s}\left(  x\right)  =e^{-sx}%
\end{equation}
gives the weighted moment
\begin{align}
m_{n}\left(  s\right)   &  =\int_{0}^{\infty}\lambda(2-q)(1-\lambda
(1-q)x)^{\frac{1}{1-q}}x^{n}e^{-sx}dx\nonumber\\
&  =\frac{\Gamma(n+1)\Gamma\left(  \frac{1}{q-1}-n-1\right)  \,_{1}%
F_{1}\left(  n+1;n+\frac{1}{1-q}+2;\frac{s}{(q-1)\lambda}\right)  }{\left[
\lambda(q-1)\right]  ^{n}\Gamma\left(  \frac{1}{q-1}-1\right)  }\nonumber\\
&  -\lambda(q-2)\left[  \lambda(q-1)\right]  ^{\frac{1}{1-q}}s^{\frac{1}%
{q-1}-1-n}\Gamma\left(  n+1+\frac{1}{1-q}\right)  \,_{1}F_{1}\left(  \frac
{1}{q-1};\frac{1}{q-1}-n;\frac{s}{(q-1)\lambda}\right)  .
\end{align}
\bigskip According to Eq. (\ref{mnrena0}), we expand the weighted moment
$m_{n}\left(  s\right)  $ at $s=0$:%
\begin{align}
m_{n}\left(  s\right)   &  =\frac{\Gamma(n+1)\Gamma\left(  -n+\frac{1}%
{q-1}-1\right)  }{\Gamma\left(  \frac{1}{q-1}-1\right)  \left[  \lambda
(q-1)\right]  ^{n}}\nonumber\\
&  +s^{\frac{1}{q-1}-1-n}\left[  \frac{(q-2)s\lambda^{\frac{1}{1-q}}%
\Gamma\left(  n+1-\frac{1}{q-1}\right)  }{(q-1)^{\frac{q}{q-1}}\left[
n\left(  q-1\right)  -1\right]  }-\lambda(q-2)\left[  \lambda(q-1)\right]
^{\frac{1}{1-q}}\Gamma\left(  \frac{1}{1-q}+n+1\right)  \right] \nonumber\\
&  +\frac{(n+1)s\Gamma(n+1)\Gamma\left(  \frac{1}{q-1}-1-n\right)  }{\left[
\lambda(q-1)\right]  ^{n}\lambda(nq-n+2q-3)\Gamma\left(  \frac{1}%
{q-1}-1\right)  }+\cdots,
\end{align}
where the zero-order term is the renormalized moment:%
\begin{equation}
m_{n}^{\text{ren}}=\frac{\Gamma(n+1)\Gamma\left(  -n+\frac{1}{q-1}-1\right)
}{\Gamma\left(  \frac{1}{q-1}-1\right)  \left[  \lambda(q-1)\right]  ^{n}}.
\end{equation}

\subsubsection{Weight function method (2)}

Taking another weight function as
\begin{equation}
g\left(  x\right)  =x^{s-1}%
\end{equation}
the weight moment is%
\begin{align}
m_{n}\left(  s\right)   &  =\int_{0}^{\infty}\lambda(2-q)\left[
1-\lambda(1-q)x\right]  ^{\frac{1}{1-q}}x^{n}x^{s-1}dx\nonumber\\
&  =\frac{\Gamma(n+s)\Gamma\left(  \frac{1}{q-1}-n-s\right)  }{\Gamma\left(
\frac{1}{q-1}-1\right)  \left[  \lambda(q-1)\right]  ^{n+s-1}}.
\end{align}
The weighted moment at $s=1$ is the renormalized $n$-th moment of the
$q$-exponential distribution:%
\begin{equation}
m_{n}^{\text{ren}}=m_{n}\left(  1\right)  =\frac{\Gamma(n+1)\Gamma\left(
\frac{1}{q-1}-n-1\right)  }{\Gamma\left(  \frac{1}{q-1}-1\right)  \left[
\lambda(q-1)\right]  ^{n}}.
\end{equation}

\subsubsection{Cut-off method}

We use the cut-off method in section \ref{IudR} to calculate the renormalized
$n$-th moment of the $q$-exponential distribution.

Cutting off the upper integration limit of the $n$-th moment gives%
\begin{equation}
m_{n}\left(  \Lambda\right)  =\int_{0}^{\Lambda}\lambda(2-q)\left[
1-\lambda(1-q)x\right]  ^{\frac{1}{1-q}}x^{n}dx. \label{qexpcut}%
\end{equation}
When $\Lambda\rightarrow\infty$, we obtain the $n$-th moment: $m_{n}%
=m_{n}\left(  \infty\right)  $. Working out the integral in Eq. (\ref{qexpcut}%
) gives%
\begin{equation}
m_{n}\left(  \Lambda\right)  =\lambda(2-q)\Lambda^{n+1}\Gamma(n+1)\,_{2}%
\tilde{F}_{1}\left(  n+1,\frac{1}{q-1};n+2;(1-q)\lambda\Lambda\right)  .
\end{equation}
In order to obtain the renormalized $n$-th moment, we expand $m_{n}\left(
\Lambda\right)  $ at $\Lambda\rightarrow\infty$:%
\begin{align}
m_{n}\left(  \Lambda\right)   &  =\frac{\Gamma(n+1)\Gamma\left(  \frac{1}%
{q-1}-1-n\right)  }{\Gamma\left(  \frac{1}{q-1}-1\right)  \left[
\lambda(q-1)\right]  ^{n}}\nonumber\\
&  -\frac{\lambda(q-2)\left[  \lambda(q-1)\right]  ^{\frac{q}{1-q}}%
\Lambda^{n+\frac{1}{1-q}}\left\{  \lambda\Lambda(q-1)-\left[  1+\frac
{q-1}{n(q-1)-1}\right]  \right\}  }{\left[  n(q-1)-1\right]  \left[
1+\frac{q-1}{n(q-1)-1}\right]  }+\cdots.
\end{align}
Dropping the divergent terms at $\Lambda\rightarrow\infty$, we obtain the
renormalized $n$-th moment:%
\begin{equation}
m_{n}^{\text{ren}}=\frac{\Gamma(n+1)\Gamma\left(  \frac{1}{q-1}-1-n\right)
}{\Gamma\left(  \frac{1}{q-1}-1\right)  \left[  \lambda(q-1)\right]  ^{n}}.
\end{equation}

\subsubsection{Characteristic function method}

We use the characteristic function method in section \ref{chm} to calculate
the renormalized $n$-th moment of the $q$-exponential distribution.

The characteristic function of the $q$-exponential distribution, by Eq.
(\textit{\ref{chfunction}}), is%
\begin{equation}
f\left(  k\right)  =-\frac{q-2}{q-1}e^{\frac{ik}{\lambda-\lambda q}}%
E_{\frac{1}{q-1}}\left(  \frac{ik}{\lambda-q\lambda}\right)  , \label{qexpC}%
\end{equation}
where $E_{n}\left(  z\right)  =\int_{1}^{\infty}\frac{e^{-zt}}{t^{n}}dt$ is
the exponential integral function\textit{. }

By the integral representation of the derivation, Eq. (\ref{mnintkdfdkn}), we
obtain%
\begin{align}
\frac{d^{n}f\left(  k\right)  }{dk^{n}}  &  =-\frac{(q-2)k^{-n-1}}{q-1}\left[
-k\Gamma\left(  1+\frac{1}{1-q}\right)  \,_{2}\tilde{F}_{2}\left(
1,1;1-n,2+\frac{1}{1-q};\frac{ik}{\lambda-q\lambda}\right)  \right.
\nonumber\\
&  \left.  -i\lambda(q-2)\Gamma\left(  \frac{1}{q-1}-1\right)  \Gamma\left(
\frac{q-2}{q-1}\right)  \left(  \frac{ik}{\lambda-\lambda q}\right)
^{\frac{1}{q-1}}\,_{1}\tilde{F}_{1}\left(  \frac{1}{q-1};\frac{1}{q-1}%
-n;\frac{ik}{\lambda-q\lambda}\right)  \right]  , \label{qexpc}%
\end{align}
where $_{p}\tilde{F}_{q}\left(  a;b;z\right)  $ is regularized hypergeometric
function, given by $_{p}F_{q}\left(  a;b;z\right)  /\left(  \Gamma
(b_{1})\cdots\Gamma(b_{q})\right)  $.

By (\ref{mndfdk}), the renormalized $n$-th moment is%
\begin{align}
m_{1}^{\text{ren}}  &  =\left(  -i\right)  \left.  \frac{df\left(  k\right)
}{dk}\right\vert _{k=0}=-\frac{1}{\lambda\left(  2q-3\right)  },\nonumber\\
m_{2}^{\text{ren}}  &  =\left(  -i\right)  ^{2}\left.  \frac{d^{2}f\left(
k\right)  }{dk^{2}}\right\vert _{k=0}=\frac{2}{\lambda^{2}(2q-3)(3q-4)}%
,\nonumber\\
&  \cdots.
\end{align}

\subsubsection{Mellin transform method: characteristic function}

We use the Mellin transform method of the characteristic function in section
\ref{MellinC} to calculate the renormalized $n$-th moment of the
$q$-exponential distribution.

Substituting the characteristic function of the $q$-exponential distribution
(\ref{qexpC}) into Eq. (\ref{powermellin}), we obtain%
\begin{align}
m_{n}^{\text{ren}}  &  =\frac{i^{n}}{\Gamma(-n)}\mathcal{M}_{-n}\left[
f\left(  k\right)  \right] \nonumber\\
&  =-\frac{\pi\csc(\pi n)}{\Gamma(-n)}\frac{\Gamma\left(  \frac{1}%
{q-1}-n-1\right)  }{\Gamma\left(  \frac{1}{q-1}-1\right)  \left[
\lambda(q-1)\right]  ^{n}}.
\end{align}
By $\Gamma(n+1)=-\pi\csc\left(  \pi n\right)  /\Gamma(-n)$, we have%
\begin{equation}
m_{n}^{\text{ren}}=\frac{\Gamma(n+1)\Gamma\left(  \frac{1}{q-1}-n-1\right)
}{\Gamma\left(  \frac{1}{q-1}-1\right)  \left[  \lambda(q-1)\right]  ^{n}}.
\end{equation}

\subsubsection{Mellin transform method: probability density function}

We use the Mellin transform method of the probability density function in
section \ref{MellinD} to calculate the renormalized $n$-th moment of the
$q$-exponential distribution.

Substituting the probability density function of the $q$-exponential
distribution, Eq. (\ref{q-exp}), into Eq. (\ref{MellinDF}), we obtain%
\begin{equation}
m_{n}^{\text{ren}}=\mathcal{M}_{n+1}\left[  p\left(  x\right)  \right]
=\frac{\Gamma(n+1)\Gamma\left(  \frac{1}{q-1}-n-1\right)  }{\Gamma\left(
\frac{1}{q-1}-1\right)  \left[  \lambda(q-1)\right]  ^{n}}.
\end{equation}

\subsubsection{Power-logarithmic moment method}

We use the power-logarithmic moment method in section \ref{LogR} to calculate
the renormalized $n$-th moment of the $q$-exponential distribution.

The power-logarithmic moment of the $q$-exponential distribution, by Eq.
(\ref{Mn-1m}), is
\begin{align}
M_{n-1}  &  =\int_{0}^{\infty}\lambda(2-q)(1-\lambda(1-q)x)^{\frac{1}{1-q}%
}x^{n-1}\ln xdx\nonumber\\
&  =\frac{\lambda(q-2)\Gamma(n)\Gamma\left(  \frac{1}{q-1}-n\right)  \left[
\psi^{(0)}\left(  \frac{1}{q-1}-n\right)  -\psi^{(0)}(n)+\ln(\lambda
(q-1))\right]  }{\Gamma\left(  \frac{1}{q-1}\right)  \left[  \lambda
(q-1)\right]  ^{n}}.
\end{align}
By the relation between the $n$-th moment $m_{n}$ and the power-logarithmic
moment $M_{n-1}$, Eq. (\ref{mnMn}), we have
\begin{equation}
m_{n}^{\text{ren}}=\frac{\Gamma(n+1)\Gamma\left(  \frac{1}{q-1}-n-1\right)
}{\Gamma\left(  \frac{1}{q-1}-1\right)  \left[  \lambda(q-1)\right]  ^{n}}.
\end{equation}

\subsubsection{Removing singularity from renormalized $n$-th moment}

Before renormalization, the $n$-th moment (\ref{qemn}) of the $q$-exponential
distribution (\ref{q-exp}) only exists when%
\begin{equation}
q<\frac{n+2}{n+1}.
\end{equation}
After renormalization, we obtain the renormalized $n$-th moment (\ref{qemnr}).
However, the renormalized moment (\ref{qemnr}) still has singularities at
\begin{align}
q  &  =1,\\
q  &  =\frac{n+2}{n+1},\text{ \ }n\geq1. \label{qsing}%
\end{align}
Here $q=1$ is the removable singularity and $q=\frac{n+2}{n+1}$ is the
essential singularity.

$q=1$ is a removable singularity. The zero-order term of the expansion at
$q=1$ of the renormalized moment (\ref{qemnr}) (i.e., the limit of
$q\rightarrow1$) is the renormalized $n$-th moment:%
\begin{equation}
m_{n}^{\text{ren}}\left(  q=1\right)  =\frac{1}{\lambda^{n}}\Gamma\left(
n+1\right)  . \label{qexpexp}%
\end{equation}
When $q=1$, the $q$-exponential distribution returns to the exponential
distribution $p\left(  x\right)  =\lambda e^{-\lambda x}$. Eq. (\ref{qexpexp})
just gives the $n$-th moment of the exponential distribution.

For the moment at the essential singularity $q=\frac{n+2}{n+1}$, we need
further renormalization. Expanding the renormalized $n$-th moment
(\ref{qemnr}) at $q=\frac{n+2}{n+1}$, we arrive at%
\begin{align}
m_{n}^{\text{ren}}  &  =-\frac{1}{\lambda^{n}}\left(  \frac{1}{1+n}\right)
^{2-n}\frac{1}{q-\frac{n+2}{n+1}}\frac{\Gamma\left(  n+1\right)  }%
{\Gamma\left(  n\right)  }\nonumber\\
&  +\frac{n}{\lambda^{n}}\left(  \frac{1}{n+1}\right)  ^{1-n}\left[  \left(
n+1\right)  \left(  \psi^{(0)}(n)+\gamma_{\text{E}}\right)  +1-n\right]
+\cdots,
\end{align}
where $\gamma_{\text{E}}=0.577216$ is Euler's constant. By the minimal
subtraction, dropping the divergent terms at $q=\frac{n+2}{n+1}$, we obtain
the renormalized $n$-th moment%
\begin{equation}
\left.  m_{n}^{\text{ren}}\right\vert _{q=\frac{n+2}{n+1}}=\frac{n}%
{\lambda^{n}}\left(  \frac{1}{n+1}\right)  ^{1-n}\left[  \left(  n+1\right)
\left(  \psi^{(0)}(n)+\gamma_{\text{E}}\right)  +1-n\right]  .
\end{equation}

\subsection{$q$-Gaussian distribution}

The probability density function of the $q$-Gaussian distribution is
\cite{domingo2017boundedness}%
\begin{equation}
p\left(  x\right)  =\left\{
\begin{array}
[c]{c}%
\frac{\sqrt{1-q}}{\sqrt{2\pi}}\left[  \frac{(q-1)(\mu-x)^{2}}{2\beta^{3-q}%
}+1\right]  ^{\frac{1}{1-q}}\frac{\Gamma\left(  \frac{1}{1-q}+\frac{3}%
{2}\right)  }{\Gamma\left(  \frac{1}{1-q}+1\right)  },\text{ }q<1,\text{
}-1<\frac{\sqrt{1-q}}{\sqrt{2}\beta}\left(  x-\mu\right)  <1,\\
\frac{1}{\sqrt{2\pi}\beta}\exp\left(  -\frac{(\mu-x)^{2}}{2\beta^{2}}\right)
,\text{ \ }q=1,\\
\frac{\sqrt{q-1}}{\sqrt{2\pi}}\left[  \frac{(q-1)(\mu-x)^{2}}{2\beta^{3-q}%
}+1\right]  ^{\frac{1}{1-q}}\frac{\Gamma\left(  \frac{1}{q-1}\right)  }%
{\Gamma\left(  \frac{1}{1-q}+\frac{1}{2}\right)  },\text{ \ }1<q<3,\\
0,\text{ \ }q>3.
\end{array}
\right.  . \label{pqgauss}%
\end{equation}

When $q=1$, the $q$-Gaussian distribution recovers the Gaussian distribution,
and the $n$-th moment of the Gaussian distribution exists. When $q=2$, the
$q$-Gaussian distribution recovers the Cauchy distribution, and though the
$n$-th moment of the Cauchy distribution does not exist, we have obtained the
renormalized $n$-th moment in section \ref{Cauchy}. Thus, what we need to
renormalize is the $n$-th moment of $1<q<3$,%
\begin{equation}
m_{n}=\int_{-\infty}^{\infty}\frac{\sqrt{q-1}}{\sqrt{2\pi}}\left[
\frac{(q-1)(\mu-x)^{2}}{2\beta^{3-q}}+1\right]  ^{\frac{1}{1-q}}\frac
{\Gamma\left(  \frac{1}{q-1}\right)  }{\Gamma\left(  \frac{1}{1-q}+\frac{1}%
{2}\right)  }x^{n}dx. \label{qgmn}%
\end{equation}

The first several moments of the $q$-Gaussian distribution, by Eq.
(\ref{pqgauss}), is
\begin{align}
m_{1}  &  =\left\{
\begin{array}
[c]{c}%
\mu,\text{ \ }1<q<2\\
\text{none, }2<q<3
\end{array}
\right.  ,\nonumber\\
m_{2}  &  =\left\{
\begin{array}
[c]{c}%
\mu^{2}+\frac{2\beta^{2}}{5-3q},\text{ \ }3q<5\\
\text{none, }\frac{5}{3}<q<3
\end{array}
\right.  ,\nonumber\\
m_{3}  &  =\left\{
\begin{array}
[c]{c}%
\mu^{3}+\frac{6\beta^{2}\mu}{5-3q},\text{ \ }2q<3\\
\text{none, }\frac{3}{2}<q<3
\end{array}
\right.  ,\nonumber\\
m_{4}  &  =\left\{
\begin{array}
[c]{c}%
\mu^{4}+\frac{12\beta^{4}}{15q^{2}-46q+35}+\frac{12\beta^{2}\mu^{2}}%
{5-3q},\text{ \ }5q<7\\
\text{none, }\frac{7}{5}<q<3
\end{array}
\right.  .
\end{align}
It can be seen that when $q>\frac{n+3}{n+1}$, the $n$-th moment does not
exist. Next, we use the renormalization method to give renormalized moments
for $q>\frac{n+3}{n+1}$.

\subsubsection{Zeta function method}

For simplicity, we only consider $\mu=0$.

The spectral zeta function of the $q$-Gaussian distribution of $1<q<3$, given
by Eqs. (\ref{pqgauss}) and (\ref{zetaint}), is%
\begin{align}
\zeta\left(  s\right)   &  =\int_{-\infty}^{\infty}\frac{\sqrt{q-1}%
\Gamma\left(  \frac{1}{q-1}\right)  \left(  \frac{q-1}{2\beta^{2}}%
x^{2}+1\right)  ^{\frac{1}{1-q}}}{\sqrt{2\pi}\beta\Gamma\left(  \frac{1}%
{q-1}-\frac{1}{2}\right)  }x^{-s}dx\nonumber\\
&  =\frac{\left[  (-1)^{-s}+1\right]  \left(  q-1\right)  ^{s/2}\Gamma\left(
\frac{1-s}{2}\right)  \Gamma\left(  \frac{1}{q-1}-\frac{1-s}{2}\right)
}{2^{s/2+1}\sqrt{\pi}\beta^{s}\Gamma\left(  \frac{1}{q-1}-\frac{1}{2}\right)
}.
\end{align}
We are concerned about the $n$-th moment. By Eq. (\ref{mnzeta}), the
renormalized $n$-th moment of the $q$-Gaussian distribution is%
\begin{align}
m_{n}^{\text{ren}}  &  =\zeta\left(  -n\right) \nonumber\\
&  =\frac{2^{n/2-1}\left[  (-1)^{n}+1\right]  \beta^{n}\Gamma\left(
\frac{n+1}{2}\right)  \Gamma\left(  \frac{1}{q-1}-\frac{1+n}{2}\right)
}{\sqrt{\pi}\left(  q-1\right)  ^{n/2}\Gamma\left(  \frac{1}{q-1}-\frac{1}%
{2}\right)  }. \label{qgmnr}%
\end{align}
Thus, the first several renormalized moments are%
\begin{align}
m_{1}^{\text{ren}}  &  =\zeta\left(  -1\right)  =0,\nonumber\\
m_{2}^{\text{ren}}  &  =\zeta\left(  -2\right)  =-\frac{2\beta^{2}}%
{5-3q},\nonumber\\
m_{3}^{\text{ren}}  &  =\zeta\left(  -3\right)  =0,\nonumber\\
m_{4}^{\text{ren}}  &  =\zeta\left(  -4\right)  =\frac{12\beta^{4}}%
{15q^{2}-46q+35}.
\end{align}

\subsubsection{Subtraction method}

The $n$-th moment of the $q$-Gaussian distribution exists at $q\leq1$, but it
does not exist at $1<q<3$ and $q>\frac{n+3}{n+1}$. Thus, we only need to deal
with the latter. We use the subtraction method in section \ref{removeinf} to
calculate the renormalized $n$-th moment. For simplicity, we only consider
$\mu=0$.

The $n$-th moment of the $q$-Gaussian distribution (\ref{qemn}) diverges at
$x\rightarrow\pm\infty$. Expand the integrand at $x\rightarrow\pm\infty$:%

\begin{align}
&  \frac{\sqrt{q-1}\Gamma\left(  \frac{1}{q-1}\right)  }{\sqrt{2\pi}%
\beta\Gamma\left(  \frac{3-q}{2(q-1)}\right)  }x^{n}x^{\frac{2}{1-q}}\left(
\frac{q-1}{2\beta^{2}}+\frac{1}{x^{2}}\right)  ^{\frac{1}{1-q}}\nonumber\\
&  =\frac{\sqrt{q-1}\Gamma\left(  \frac{1}{q-1}\right)  }{\sqrt{2\pi}%
\beta\Gamma\left(  \frac{3-q}{2(q-1)}\right)  }\sum_{j=0}^{\infty}\left(
\begin{array}
[c]{c}%
\frac{1}{1-q}\\
j
\end{array}
\right)  \left(  \frac{q-1}{2\beta^{2}}\right)  ^{\frac{1}{1-q}-j}\left(
\frac{1}{x^{2}}\right)  ^{j+\frac{1}{q-1}}\frac{1}{x^{n}},
\end{align}
where $\left(
\begin{array}
[c]{c}%
m\\
n
\end{array}
\right)  =\frac{\Gamma\left(  m+1\right)  }{\Gamma\left(  n+1\right)
\Gamma\left(  m-n+1\right)  }$\ is the binomial expansion coefficient.
Following section \ref{removeinf}, the divergence is caused by the term
$2j-n+\frac{2}{q-1}\leq1$, i.e.,%
\begin{equation}
j\leq\frac{1}{2}+\frac{n}{2}-\frac{1}{q-1}.
\end{equation}
By Eq. (\ref{mreninf}), we obtain the renormalized $n$-th moment%
\begin{align}
m_{n}^{\text{ren}}  &  =\int_{-1}^{1}f_{q}\left(  x\right)  dx+\int_{-\infty
}^{-1}\left[  f_{q}\left(  x\right)  -\sum_{j=0}^{j\leq\frac{1}{2}+\frac{n}%
{2}-\frac{1}{q-1}}c_{j}\left(  \frac{1}{x^{2}}\right)  ^{j+\frac{1}{q-1}}%
\frac{1}{x^{n}}\right]  dx\nonumber\\
&  +\sum_{j=0}^{j\leq\frac{1}{2}+\frac{n}{2}-\frac{1}{q-1}}\frac{\left(
-1\right)  ^{-n-1}c_{j}}{1-\left(  2j-n+\frac{2}{q-1}\right)  }+i\pi
c_{\text{$-1$}}+\int_{1}^{\infty}\left[  f_{q}\left(  x\right)  -\sum
_{j=0}^{j\leq\frac{1}{2}+\frac{n}{2}-\frac{1}{q-1}}c_{j}\left(  \frac{1}%
{x^{2}}\right)  ^{j+\frac{1}{q-1}}\frac{1}{x^{n}}\right]  dx\\
&  -\sum_{j=0}^{j\leq\frac{1}{2}+\frac{n}{2}-\frac{1}{q-1}}\frac{c_{j}%
}{1-\left(  2j-n+\frac{2}{q-1}\right)  },
\end{align}
where%
\begin{align}
f_{q}\left(  x\right)   &  =\frac{\sqrt{q-1}\Gamma\left(  \frac{1}%
{q-1}\right)  \left[  \frac{(q-1)x^{2}}{2\beta^{2}}+1\right]  ^{\frac{1}{1-q}%
}}{\sqrt{2\pi}\beta\Gamma\left(  \frac{3-q}{2(q-1)}\right)  }x^{n},\nonumber\\
c_{j}  &  =\frac{\sqrt{q-1}\Gamma\left(  \frac{1}{q-1}\right)  }{\sqrt{2\pi
}\beta\Gamma\left(  \frac{3-q}{2(q-1)}\right)  }\left(
\begin{array}
[c]{c}%
\frac{1}{1-q}\\
j
\end{array}
\right)  \left(  \frac{q-1}{2\beta^{2}}\right)  ^{\frac{1}{1-q}-j}.
\end{align}
Working out the integral gives%
\begin{align}
m_{n}^{\text{ren}}  &  =\frac{2}{\sqrt{2\pi}}\left[  \left(  -1\right)
^{n}+1\right]  \sqrt{q-1}\left[  \Gamma\left(  \frac{1}{q-1}\right)
\frac{\text{ }_{2}F_{1}\left(  \frac{n+1}{2},\frac{1}{q-1},\frac{n+3}{2}%
,\frac{1-q}{2\beta^{2}}\right)  }{\beta\left(  n+1\right)  \Gamma\left(
\frac{3-q}{2(q-1)}\right)  }\right. \nonumber\\
&  \left.  -\frac{\left(  \frac{q-1}{2\beta^{2}}\right)  ^{\frac{q-2}{q-1}%
}\beta\text{ }_{2}F_{1}\left(  -\frac{1}{2}-\frac{n}{2}+\frac{1}{q-1},\frac
{1}{q-1},\frac{1}{2}-\frac{n}{2}+\frac{1}{q-1},\frac{2\beta^{2}}{1-q}\right)
}{q-3+n\left(  q-1\right)  }\right]  . \label{qguass}%
\end{align}
It can be verified that this renormalized $n$-th moment is the same as that
given by the zeta function method.

From Eq. (\ref{qguass}), for example, the first and second moments are%
\begin{align}
m_{1}^{\text{ren}}  &  =0,\nonumber\\
m_{2}^{\text{ren}}  &  =-\frac{2\beta^{2}}{5-3q}.
\end{align}

\subsubsection{Cut-off method}

We use the cut-off method in section \ref{IudR} to calculate the renormalized
$n$-th moment of the $q$-Gaussian distribution.

Cutting off the upper integration limit of the $n$-th moment gives%
\begin{equation}
m_{n}\left(  \Lambda\right)  =\int_{-\Lambda}^{\Lambda}\frac{\sqrt{q-1}}%
{\sqrt{2\pi}}\left[  \frac{(q-1)}{2\beta^{3-q}}x^{2}+1\right]  ^{\frac{1}%
{1-q}}\frac{\Gamma\left(  \frac{1}{q-1}\right)  }{\Gamma\left(  \frac{1}%
{1-q}+\frac{1}{2}\right)  }x^{n}dx. \label{qgausecut}%
\end{equation}
When $\Lambda\rightarrow\infty$, we obtain the $n$-th moment: $m_{n}%
=m_{n}\left(  \infty\right)  $. Working out the integral in Eq.
(\ref{qgausecut}) gives
\begin{equation}
m_{n}\left(  \Lambda\right)  =\frac{\Lambda\sqrt{q-1}\left(  \Lambda
^{n}+(-\Lambda)^{n}\right)  \Gamma\left(  \frac{1}{q-1}\right)  \,_{2}%
F_{1}\left(  \frac{n+1}{2},\frac{1}{q-1};\frac{n+3}{2};-\frac{(q-1)\Lambda
^{2}}{2\beta^{2}}\right)  }{\sqrt{2\pi}\beta(n+1)\Gamma\left(  \frac{1}%
{q-1}-\frac{1}{2}\right)  }.
\end{equation}
In order to obtain the renormalized $n$-th moment, we expand $m_{n}\left(
\Lambda\right)  $ at $\Lambda\rightarrow\infty$:%
\begin{align}
m_{n}\left(  \Lambda\right)   &  =\frac{2^{n/2-1}\left[  (-1)^{n}+1\right]
\beta^{n}\Gamma\left(  \frac{n+1}{2}\right)  \Gamma\left(  \frac{1}{q-1}%
-\frac{1+n}{2}\right)  }{\sqrt{\pi}\Gamma\left(  \frac{1}{q-1}-\frac{1}%
{2}\right)  \left(  q-1\right)  ^{n/2}}\nonumber\\
&  +\Lambda^{n+1-\frac{2}{q-1}}\frac{2^{\frac{1}{q-1}+\frac{1}{2}}\left[
1+(-1)^{n}\right]  \Gamma\left(  \frac{1}{q-1}\right)  \beta^{\frac{2}{q-1}%
-1}\left\{  (q-1)^{3}(n-1)-\frac{2\beta^{2}}{\Lambda^{2}}\left[
n(q-1)+q-3\right]  \right\}  }{(q-1)^{\frac{3q-1}{2\left(  q-1\right)  }%
}(n-1)\left[  n(q-1)+q-3\right]  2\sqrt{\pi}\Gamma\left(  \frac{1}{q-1}%
-\frac{1}{2}\right)  }+\cdots.
\end{align}
Dropping the divergent terms at $\Lambda\rightarrow\infty$, we obtain the
renormalized $n$-th moment,%
\begin{equation}
m_{n}^{\text{ren}}=\frac{2^{n/2-1}\left(  (-1)^{n}+1\right)  \beta^{n}%
\Gamma\left(  \frac{n+1}{2}\right)  \Gamma\left(  \frac{1}{q-1}-\frac{1+n}%
{2}\right)  }{\sqrt{\pi}\Gamma\left(  \frac{1}{q-1}-\frac{1}{2}\right)
\left(  q-1\right)  ^{n/2}}.
\end{equation}

\subsubsection{Characteristic function method}

We use the characteristic function method in section \ref{chm} to calculate
the renormalized $n$-th moment of the $q$-Gaussian distribution.

The characteristic function of the $q$-Gaussian distribution, by Eq.
(\ref{chfunction}), is%
\begin{align}
f\left(  k\right)   &  =\text{ }_{0}F_{1}\left(  ;\frac{3}{2}-\frac{1}%
{q-1};\frac{k^{2}\beta^{2}}{2(q-1)}\right)  -\frac{\beta^{\frac{3-q}{q-1}%
}\left(  q-1\right)  ^{\frac{3q-5}{2\left(  q-1\right)  }}}{2^{\frac{q-2}%
{q-1}}\sqrt{2\pi}}\Gamma\left(  \frac{2\left(  q-2\right)  }{q-1}\right)
\nonumber\\
&  \times\Gamma\left(  \frac{q}{q-1}\right)  \left(  i^{\frac{q-3}{q-1}%
}-i^{\frac{q+1}{q-1}}\right)  \,_{0}\tilde{F}_{1}\left(  ;\frac{1}{2}+\frac
{1}{q-1};\frac{k^{2}\beta^{2}}{2(q-1)}\right)  k^{\frac{3-q}{q-1}}.
\label{q-GaussianC}%
\end{align}
By the integral representation of the derivation, Eq. (\ref{mnintkdfdkn}), we
obtain%
\begin{align}
\frac{d^{n}f\left(  k\right)  }{dk^{n}}  &  =2^{n}\sqrt{\pi}\,_{2}\tilde
{F}_{3}\left(  \frac{1}{2},1;\frac{1-n}{2},1-\frac{n}{2},\frac{3}{2}+\frac
{1}{1-q};\frac{k^{2}\beta^{2}}{2(q-1)}\right)  \Gamma\left(  \frac{3}{2}%
+\frac{1}{1-q}\right)  k^{-n}\nonumber\\
&  -\frac{2^{n}\beta^{\frac{3-q}{q-1}}}{\sqrt{\pi}\left(  \frac{q-1}%
{2}\right)  ^{\frac{5-3q}{2\left(  q-1\right)  }}}\left(  i^{\frac{q-3}{q-1}%
}-i^{\frac{q+1}{q-1}}\right)  \Gamma\left(  \frac{q-2}{\left(  q-1\right)
/2}\right)  \Gamma\left(  \frac{q}{q-1}\right) \nonumber\\
&  \times\,_{1}\tilde{F}_{2}\left(  \frac{1}{q-1};\frac{1}{q-1}-\frac{n}%
{2},\frac{1}{q-1}+\frac{1-n}{2};\frac{k^{2}\beta^{2}}{2(q-1)}\right)
k^{\frac{3-q}{q-1}-n}. \label{qGaussc}%
\end{align}

$q=\frac{3}{2}$ and $q=\frac{4}{3}$ are not singularities, and the moment can
be obtained directly from Eq. (\ref{qGaussc}): at $q=\frac{3}{2}$%
\begin{equation}
m_{1}^{\text{ren}}=\left(  -i\right)  \left.  \frac{df\left(  k\right)  }%
{dk}\right\vert _{k=0}=0;
\end{equation}
at $q=\frac{4}{3}$%
\begin{equation}
m_{2}^{\text{ren}}=\left(  -i\right)  ^{2}\left.  \frac{d^{2}f\left(
k\right)  }{dk^{2}}\right\vert _{k=0}=-2\beta^{2}.
\end{equation}
This means that $q=\frac{3}{2}$ and $q=\frac{4}{3}$ are not true
singularities, and the renormalized $n$-th moment can be obtained by
analytical continuation directly.

$q=2$\ and $q=\frac{5}{3}$ are singularities, which need to be renormalized.

For $q=2$,%
\begin{equation}
m_{1}^{\text{ren}}=\left(  -i\right)  \left.  \frac{df\left(  k\right)  }%
{dk}\right\vert _{k=0}=i\sqrt{2}\beta;
\end{equation}
for $q=\frac{5}{3}$,
\begin{equation}
m_{2}^{\text{ren}}=\left(  -i\right)  ^{2}\left.  \frac{d^{2}f\left(
k\right)  }{dk^{2}}\right\vert _{k=0}=3\beta^{2}\left(  \ln\frac{2}{\sqrt
{3}\beta}-1-i\frac{\pi}{2}\right)  .
\end{equation}
This means that $q=2$ and $q=\frac{5}{3}$ are true singularities, and the
finite moments cannot be obtained through analytical continuation.

\subsubsection{Mellin transform method: characteristic function}

We use the Mellin transform method of the characteristic function in section
\ref{MellinC} to calculate the renormalized $n$-th moment of the $q$-Gaussian distribution.

Substituting the characteristic function of the $q$-Gaussian distribution
(\ref{q-GaussianC}) into Eq. (\ref{powermellin}) gives immediately%
\begin{align}
m_{n}^{\text{ren}}  &  =\frac{i^{n}}{\Gamma(-n)}\mathcal{M}_{-n}\left[
f\left(  k\right)  \right] \nonumber\\
&  =\frac{2^{n/2-1}\left(  (-1)^{n}+1\right)  \beta^{n}\Gamma\left(
\frac{n+1}{2}\right)  \Gamma\left(  \frac{1}{q-1}-\frac{1+n}{2}\right)
}{\sqrt{\pi}\Gamma\left(  \frac{1}{q-1}-\frac{1}{2}\right)  \left(
q-1\right)  ^{n/2}}.
\end{align}

\subsubsection{Mellin transform method: probability density function}

We use the Mellin transform method of the probability density function in
section \ref{MellinD} to calculate the renormalized $n$-th moment of the
$q$-Gaussian distribution.

Substituting the probability density function of the $q$-Gaussian distribution
(\ref{pqgauss}) into Eq. (\ref{MellinDF}) gives%
\begin{align}
m_{n}^{\text{ren}}  &  =\left(  -1\right)  ^{n}\mathcal{M}_{n+1}\left[
p\left(  -x\right)  \right]  +\mathcal{M}_{n+1}\left[  p\left(  x\right)
\right] \nonumber\\
&  =\frac{2^{n/2-1}\left(  (-1)^{n}+1\right)  \beta^{n}\Gamma\left(
\frac{n+1}{2}\right)  \Gamma\left(  \frac{1}{q-1}-\frac{1+n}{2}\right)
}{\sqrt{\pi}\Gamma\left(  \frac{1}{q-1}-\frac{1}{2}\right)  \left(
q-1\right)  ^{n/2}}.
\end{align}

\subsubsection{Removing singularity from renormalized $n$-th moment}

Before renormalization, the $n$-th moment (\ref{qemn}) of the $q$-exponential
distribution (\ref{pqgauss}) only exists when%
\begin{equation}
q<\frac{n+3}{n+1}.
\end{equation}
After renormalization, we obtain the renormalized $n$-th moment
(\ref{Cauchymnr}). However, the renormalized moment (\ref{qgmnr}) still has
singularities at%
\begin{align}
q  &  =1,\\
q  &  =\frac{n+3}{n+1},\text{ \ \ }n\geq1.
\end{align}
Here $q=1$ is the removable singularity and $q=\frac{n+3}{n+1}$ is the
essential singularity.

$q=1$ is a removable singularity, the zero-order term of the expansion at
$q=1$ of the renormalized moment (\ref{qgmnr}) (i.e., the limit of
$q\rightarrow1$) is the renormalized $n$-th moment:%
\begin{equation}
\left.  m_{n}^{\text{ren}}\right\vert _{q=1}=\frac{2^{\frac{n}{2}-1}}%
{\sqrt{\pi}}\left[  \left(  -1\right)  ^{n}+1\right]  \beta^{n}\Gamma\left(
\frac{1}{2}+\frac{n}{2}\right)  . \label{qGaussq1}%
\end{equation}
\bigskip When $q=1$, the $q$-exponential distribution returns to the
exponential distribution $p\left(  x\right)  =\frac{1}{\sqrt{2\pi}}%
e^{-x^{2}/2}$. Eq. (\ref{qGaussq1}) just gives the $n$-th moment of the normal distribution.

For the moment at the essential singularity $q=\frac{n+3}{n+1}$, we need
further renormalization. Expanding the renormalized $n$-th moment
(\ref{qgmnr}) at $q=\frac{n+3}{n+1}$, we arrive at%
\begin{align}
m_{n}^{\text{ren}}  &  =-\frac{2\left[  \left(  -1\right)  ^{n}+1\right]
\left(  \frac{1}{n+1}\right)  ^{2-\frac{n}{2}}\beta^{n}}{\sqrt{\pi}\left(
q-\frac{n+3}{n+1}\right)  }\frac{\Gamma\left(  \frac{n+1}{2}\right)  }%
{\Gamma\left(  \frac{n}{2}\right)  }\nonumber\\
&  -\frac{\left[  \left(  -1\right)  ^{n}+1\right]  \left(  \frac{1}%
{n+1}\right)  ^{1-\frac{n}{2}}\beta^{n}\Gamma\left(  \frac{n+1}{2}\right)
\left[  (n+1)\left(  \psi^{(0)}\left(  \frac{n}{2}\right)  +\gamma_{\text{E}%
}\right)  +2-n\right]  }{2\sqrt{\pi}\Gamma\left(  \frac{n}{2}\right)  }%
+\cdots.
\end{align}
By the minimal subtraction, dropping the divergent terms at $q=\frac{n+3}%
{n+1}$, we obtain the renormalized $n$-th moment:%
\begin{equation}
\left.  m_{n}^{\text{ren}}\right\vert _{q=\frac{n+3}{n+1}}=-\frac{\left[
\left(  -1\right)  ^{n}+1\right]  \left(  \frac{1}{n+1}\right)  ^{1-\frac
{n}{2}}\beta^{n}\Gamma\left(  \frac{n+1}{2}\right)  \left[  (n+1)\left(
\psi^{(0)}\left(  \frac{n}{2}\right)  +\gamma_{\text{E}}\right)  +2-n\right]
}{2\sqrt{\pi}\Gamma\left(  \frac{n}{2}\right)  }.
\end{equation}

It can be directly verified that this result is consistent with other
renormalization methods.

\section{Nonpositive integer power moment \label{noninteger}}

The $n$-th moment is the positive integer power moment. The renormalization
schemes are essentially based on analytical continuation. Moreover, the
analytical continuation can give not only positive integer power moments but
also real and even complex moments. For example, the complex moments of the
Cauchy distribution, Levy distribution, $q$-exponential distribution, and
$q$-Gaussian distribution are%
\begin{align}
m_{z}^{\text{Cauchy}}  &  =e^{iz\pi/2},\nonumber\\
m_{z}^{\text{Levy}}  &  =\frac{1}{2^{z}\sqrt{\pi}}\Gamma\left(  \frac{1}%
{2}-z\right)  ,\nonumber\\
m_{z}^{q\text{-exp}}  &  =\frac{1}{\left[  \lambda(q-1)\right]  ^{z}}%
\frac{\Gamma(z+1)\Gamma\left(  \frac{1}{q-1}-1-z\right)  }{\Gamma\left(
\frac{1}{q-1}-1\right)  },\nonumber\\
m_{z}^{q\text{-Gauss}}  &  =\frac{2^{z/2-1}\left(  1+e^{i\pi z}\right)
\Gamma\left(  \frac{z+1}{2}\right)  \beta^{z}\Gamma\left(  \frac{1}{q-1}%
-\frac{1+z}{2}\right)  }{\sqrt{\pi}\Gamma\left(  \frac{1}{q-1}-\frac{1}%
{2}\right)  \left(  q-1\right)  ^{z/2}}.
\end{align}
The renormalization treatment, besides giving moments for the distributions
that have no moment, also gives the non-positive integer moment of the
distribution that has the positive integer power moment.

\textit{Normal distribution} $p\left(  x\right)  =\frac{1}{\sqrt{2\pi}%
}e^{-\frac{x^{2}}{2}}$. The standard normal distribution does not have
negative $n$-th moments. However, the complex moment after analytical
continuation is%

\begin{equation}
m_{z}^{\text{Normal}}=\frac{2^{z/2-1}}{\sqrt{\pi}}\left(  1+e^{i\pi z}\right)
\Gamma\left(  \frac{z+1}{2}\right)  . \label{mznormal}%
\end{equation}
$m_{z}^{\text{Normal}}$ is ill-defined at negative integers. In order to
obtain the negative $n$-th moment, we use the renormalization treatment.
Expanding the complex moment (\ref{mznormal}) at $z=-1$, $-2$, $-3$, and $4$,
respectively, and dropping the divergent term, we obtain the following
negative $n$-th moment:%
\begin{equation}
m_{-1}^{\text{Normal}}=-i\sqrt{\frac{\pi}{2}},\text{ }m_{-2}^{\text{Normal}%
}=-1,\text{ }m_{-3}^{\text{Normal}}=\frac{i}{2}\sqrt{\frac{\pi}{2}},\text{
}m_{-4}^{\text{Normal}}=\frac{1}{3},\ldots.
\end{equation}

\textit{Student's t-distribution} $p\left(  x\right)  =\frac{\left(  \frac
{\nu}{\nu+x^{2}}\right)  ^{\frac{\nu+1}{2}}}{\sqrt{\nu}B\left(  \frac{\nu}%
{2},\frac{1}{2}\right)  }$. The complex moment after analytical continuation
is%
\begin{equation}
m_{z}=\frac{(-1)^{z}\left(  1+e^{-i\pi z}\right)  \nu^{z/2}\Gamma\left(
\frac{1+z}{2}\right)  \Gamma\left(  \frac{\nu-z}{2}\right)  }{2\sqrt{\pi
}\Gamma\left(  \frac{\nu}{2}\right)  }. \label{mzt}%
\end{equation}
This result still can not give the negative $n$-th moment. Expanding the
complex moment (\ref{mzt}) at $z=-1$, $-2$, $-3$, and $4$, respectively, and
dropping the divergence term, we obtain the negative $n$-th moment:
\begin{equation}
m_{-1}^{\text{T}}=-\frac{i\pi}{\sqrt{\nu}B\left(  \frac{\nu}{2},\frac{1}%
{2}\right)  },\text{ }m_{-2}^{\text{T}}=-1,\text{ }m_{-3}^{\text{T}}%
=\frac{i\sqrt{\pi}\Gamma\left(  \frac{\nu+3}{2}\right)  }{\nu^{3/2}%
\Gamma\left(  \frac{\nu}{2}\right)  },\text{ }m_{-4}^{\text{T}}=\frac{2+\nu
}{3\nu},\ldots.
\end{equation}

\textit{Laplace distribution} $p\left(  x\right)  =\frac{\lambda}%
{2}e^{-\lambda\left\vert x-\mu\right\vert }$. For $\mu=0$, the complex order
moment after analytical continuation is
\begin{equation}
m_{z}=\frac{\lambda e^{i\pi z}+1}{2\lambda^{z}}\Gamma(1+z).
\end{equation}
This result can not give the negative $n$-th moment of the Laplace
distribution. Expanding the complex moment (\ref{mzt}) at $z=-1$, $-2$, $-3$,
and $4$, respectively, and dropping the divergent term, we obtain the negative
$n$-th moment,
\begin{align}
m_{-1}^{\text{Laplace}}  &  =-i\frac{\pi}{2}\lambda,\text{ }m_{-2}%
^{\text{Laplace}}=\lambda^{2}\left(  \ln\lambda-i\frac{\pi}{2}+\gamma
_{\text{E}}-1\right)  ,\nonumber\\
m_{-3}^{\text{Laplace}}  &  =-i\frac{\pi}{4}\lambda^{3},\text{ }%
m_{-4}^{\text{Laplace}}=\frac{\lambda^{4}}{6}\left(  \ln\lambda-i\frac{\pi}%
{2}+\gamma_{\text{E}}-\frac{11}{6}\right)  \ldots.
\end{align}

\section{Calculating logarithmic moment from power moment \label{mLog}}

In this section, we provide a method to calculate the logarithmic moment
(\ref{logM}) from the power moment. However, for some distributions, the
logarithmic moment exists, but the power moment does not. In such cases, we
can use the renormalized power moment given in this paper to calculate the
logarithmic moment. In particular, this method can directly verify the
validity of the renormalized power moment by comparing the logarithmic moments
calculated by definition (\ref{logM}) with the logarithmic moments calculated
by the renormalized power moment.

\subsection{Power moment scheme}

If a distribution has power moments, it also has logarithmic moments. If the
power moment does not exist, we can first obtain the renormalized power moment
by the renormalization procedure in section \ref{powerM} and then calculate
the logarithmic moment from the renormalized power moment. Moreover,
calculating the logarithmic moment from the power moment provides an approach
for calculating the logarithmic moment.

The second characteristic function of distribution of the probability density
function $p\left(  x\right)  $,%
\begin{equation}
\phi\left(  \sigma\right)  =\int_{-\infty}^{\infty}p\left(  x\right)
x^{\sigma-1}dx, \label{IIch}%
\end{equation}
is a moment generating function of logarithmic moment
\cite{nicolas2002introduction}. The power moment can be expressed as
\begin{equation}
m_{n}=\int_{-\infty}^{\infty}p\left(  x\right)  x^{n}dx=\phi\left(
n+1\right)  . \label{mnphi}%
\end{equation}

By the relation (\ref{Dlog}) and
\begin{equation}
\left.  \frac{d^{n}x^{\sigma}}{d\sigma^{n}}\right\vert _{\sigma=0}=\left.
x^{\sigma}\ln^{n}x\right\vert _{\sigma=0}=\ln^{n}x,
\end{equation}
we have%
\begin{align}
\widetilde{m}_{n}  &  =\int_{-\infty}^{\infty}p\left(  x\right)  \ln
^{n}xdx\nonumber\\
&  =\int_{-\infty}^{\infty}p\left(  x\right)  \left.  \frac{d^{n}x^{\sigma-1}%
}{d\sigma^{n}}\right\vert _{\sigma-1=0}dx=\left.  \frac{d^{n}}{d\sigma^{n}%
}\int_{-\infty}^{\infty}p\left(  x\right)  x^{\sigma-1}dx\right\vert
_{\sigma=1}\nonumber\\
&  =\left.  \frac{d^{n}\phi\left(  \sigma\right)  }{d\sigma^{n}}\right\vert
_{\sigma=1}.
\end{align}
Compared with Eq. (\ref{mnphi}), we have%
\begin{equation}
\widetilde{m}_{n}=\left.  \frac{d^{n}m_{\sigma-1}}{d\sigma^{n}}\right\vert
_{\sigma=1}. \label{logmn}%
\end{equation}
By the relation (\ref{logmn}), the logarithmic moment can be obtained from the
renormalized power moment:%
\begin{equation}
\widetilde{m}_{n}^{\text{ren}}=\left.  \frac{d^{n}m_{\sigma-1}^{\text{ren}}%
}{d\sigma^{n}}\right\vert _{\sigma=1}. \label{logpower}%
\end{equation}

\subsection{Zeta function method}

Observing the expression of the zeta function, Eq. (\ref{zetaint}), we can see
that the second characteristic function (\ref{IIch}) can be expressed by the
zeta function:%
\begin{equation}
\phi\left(  \sigma\right)  =\zeta\left(  1-\sigma\right)  .
\end{equation}
By Eq. (\ref{logmn}), the logarithmic moment can be expressed as the zeta
function:%
\begin{equation}
\widetilde{m}_{n}=\left.  \frac{d^{n}\zeta\left(  1-\sigma\right)  }%
{d\sigma^{n}}\right\vert _{\sigma=1}.
\end{equation}
In this way, the zeta function renormalization scheme can be directly applied
to calculate the logarithmic moment.

This result shows that the zeta function ($n$ is substituted by $-n$) is the
moment generating function of the logarithmic moment.

\subsection{Characteristic function method}

The first characteristic function is the moment generating function of power
moments, and the second characteristic function is the moment generating
function of logarithmic moment \cite{nicolas2002introduction}. In this
section, we present a method of obtaining the logarithmic moment from the
first characteristic function.

By Eqs. (\ref{mndfdk}) and (\ref{mndfdk}), we directly obtain the renormalized
logarithmic moment,
\begin{align}
\widetilde{m}_{n}^{\text{ren}}  &  =\left.  \frac{d^{n}m_{\sigma
-1}^{\text{ren}}}{d\sigma^{n}}\right\vert _{\sigma=1}\nonumber\\
&  =\left.  \frac{d^{n}}{d\sigma^{n}}\left[  \left(  -i\right)  ^{\sigma
-1}\left.  \frac{d^{\sigma-1}}{dk^{\sigma-1}}f\left(  k\right)  \right\vert
_{k=0}\right]  \right\vert _{\sigma=1}. \label{powerlog}%
\end{align}
The logarithmic moment can be obtained from the first characteristic function.
It can be seen that $\left(  -i\right)  ^{\sigma-1}\left.  \frac{d^{\sigma
-1}f\left(  k\right)  }{dk^{\sigma-1}}\right\vert _{k=0}=\phi\left(
\sigma\right)  $ in Eq. (\ref{powerlog}) is the second characteristic
function, which is also the generating function of logarithmic moment.

\section{Calculating logarithmic moment from power moment: example
\label{mLogEx}}

In section \ref{mLog}, we provide a method for calculating logarithmic moments
from power moments. Below we provide some examples.

\subsection{Cauchy distribution}

In principle, the logarithmic moment of the Cauchy distribution can be
directly calculated by definition (\ref{logM}), but here we calculate the
logarithmic moment from the renormalized power moment. The renormalized power
moment of the Cauchy distribution is given by Eq. (\ref{Cauchymnr}).

By the relation between power moment and logarithmic moment, Eq.
(\ref{powerlog}), we can calculate the logarithmic moment from the $n$-th
power moment:%
\begin{align}
\widetilde{m}_{n}  &  =\left.  \frac{d^{n}m_{\sigma-1}}{d\sigma^{n}%
}\right\vert _{\sigma=1}=\left(  \frac{\pi}{2}\right)  ^{n}e^{in\pi
/2}\nonumber\\
&  =\left(  \frac{i\pi}{2}\right)  ^{n},
\end{align}
which is consistent with the result obtained from the definition of the
logarithmic moment, Eq. (\ref{logM}).

The first several logarithmic moments then read $\widetilde{m}_{1}=\frac{i}%
{2}\pi$, $\widetilde{m}_{2}=-\frac{\pi^{2}}{4}$, $\widetilde{m}_{3}=-\frac
{i}{8}\pi^{3}$, and $\widetilde{m}_{4}=\frac{\pi^{4}}{16}$.

\subsection{Levy distribution}

The renormalized power moment of the Levy distribution (\ref{Levypdf}) is
given by Eq. (\ref{Levymnr}).

By the relation between $n$-th power moments and logarithmic moments, Eq.
(\ref{powerlog}), we obtain the first-order logarithmic moment:%
\begin{equation}
\widetilde{m}_{1}=\left.  \frac{dm_{\sigma-1}}{d\sigma}\right\vert _{\sigma
=1}=\gamma_{\text{E}}+\ln2,
\end{equation}
the second-order logarithmic moment%
\begin{equation}
\widetilde{m}_{2}=\left.  \frac{d^{2}m_{\sigma-1}}{d\sigma^{2}}\right\vert
_{\sigma=1}=\frac{\pi^{2}}{2}+\ln^{2}2+\gamma_{\text{E}}\left(  \gamma
_{\text{E}}+2\ln2\right)  ,
\end{equation}
the third-order logarithmic moment%
\begin{equation}
\widetilde{m}_{3}=\left.  \frac{d^{3}m_{\sigma-1}}{d\sigma^{3}}\right\vert
_{\sigma=1}=14\zeta(3)+\left(  \gamma_{\text{E}}+\ln2\right)  ^{3}+\frac{3}%
{2}\pi^{2}\left(  \gamma_{\text{E}}+\ln2\right)  ,
\end{equation}
and the fourth logarithmic moment%
\begin{equation}
\widetilde{m}_{4}=\left.  \frac{d^{4}m_{\sigma-1}}{d\sigma^{4}}\right\vert
_{\sigma=1}=\left(  \gamma_{\text{E}}+\ln2\right)  ^{4}+3\pi^{2}%
(\gamma_{\text{E}}+\ln2)^{2}+56\zeta(3)(\gamma_{\text{E}}+\ln2)+\frac{7\pi
^{4}}{4},
\end{equation}
where $\zeta(s)$ is the Riemann zeta function. $\allowbreak$

These results are consistent with the result obtained by definition
(\ref{logM}).

\subsection{$q$-exponential distribution}

The renormalized $n$-th power moment of the $q$-exponential distribution
(\ref{q-exp}) for $q>1$ and $x\geq0$ is given by\ Eq. (\ref{qemnr}).

By the relation between the $n$-th power moment and the logarithmic moment
(\ref{powerlog}), we obtain the first-order logarithmic moment:
\begin{equation}
\widetilde{m}_{1}=\left.  \frac{dm_{\sigma-1}}{d\sigma}\right\vert _{\sigma
=1}=-H_{\frac{1}{q-1}-2}-\ln(\lambda(q-1))
\end{equation}
with $H_{n}=\sum_{i=1}^{n}\frac{1}{i}$ the $n$-th harmonic number
\cite{olver2010nist}, the second-order logarithmic moment%
\begin{equation}
\widetilde{m}_{2}=\left.  \frac{d^{2}m_{\sigma-1}}{d\sigma^{2}}\right\vert
_{\sigma=1}=\left[  H_{\frac{1}{q-1}-2}+\ln(\lambda(q-1))\right]  ^{2}%
+\psi^{(1)}\left(  \frac{1}{q-1}-1\right)  +\frac{\pi^{2}}{6}%
\end{equation}
with $\psi^{(n)}\left(  z\right)  =\frac{d^{n}\psi\left(  z\right)  }{dz^{n}}$
the $n$-th derivative of the digamma function, $\psi\left(  z\right)
=\frac{\Gamma^{\prime}\left(  z\right)  }{\Gamma\left(  z\right)  }$ the
logarithmic derivative of the gamma function, and the third logarithmic
moment
\begin{align}
\widetilde{m}_{3}  &  =\left.  \frac{d^{3}m_{\sigma-1}}{d\sigma^{3}%
}\right\vert _{\sigma=1}\nonumber\\
=  &  -\frac{1}{2}\left[  \ln(\lambda(q-1))+\gamma_{\text{E}}\right]  \left\{
2\ln(\lambda(q-1))\left[  \ln(\lambda(q-1))+2\gamma_{\text{E}}\right]
+\pi^{2}+2\gamma_{\text{E}}^{2}\right\}  -2\zeta(3)\nonumber\\
&  -\frac{1}{2}\psi^{(0)}\left(  \frac{1}{q-1}-1\right)  \left\{  6\left[
\ln(\lambda(q-1))+2\gamma_{\text{E}}\right]  \ln(\lambda(q-1))+\pi^{2}%
+6\gamma_{\text{E}}^{2}\right\} \nonumber\\
&  -3\psi^{(0)}\left(  \frac{1}{q-1}-1\right)  ^{2}\left[  \ln(\lambda
(q-1))+\gamma_{\text{E}}\right]  -\psi^{(0)}\left(  \frac{1}{q-1}-1\right)
^{3}\nonumber\\
&  -3\psi^{(1)}\left(  \frac{1}{q-1}-1\right)  \left[  H_{\frac{1}{q-1}-2}%
+\ln(\lambda(q-1))\right]  -\psi^{(2)}\left(  \frac{1}{q-1}-1\right)  .
\end{align}
The above result is consistent with the result obtained from the definition
(\ref{logM}).

\subsection{$q$-Gaussian distribution}

The renormalized $n$-th power moment of the $q$-Gaussian distribution
(\ref{pqgauss}) is Eq. (\ref{qgmnr}), which is valid for all values of $q$.

By the relation between $n$-th power moment and logarithmic moment
(\ref{powerlog}), we obtain the first-order logarithmic moment:%
\begin{align}
\widetilde{m}_{1}  &  =\left.  \frac{dm_{\sigma-1}}{d\sigma}\right\vert
_{\sigma=1}\nonumber\\
&  =-\frac{1}{2}\left[  H_{\frac{1}{q-1}-\frac{3}{2}}+\ln\frac{2(q-1)}%
{\beta^{2}}-i\pi\right]  ,
\end{align}
the second-order logarithmic moment%
\begin{align}
\widetilde{m}_{2}  &  =\left.  \frac{d^{2}m_{\sigma-1}}{d\sigma^{2}%
}\right\vert _{\sigma=1}\nonumber\\
&  =\frac{1}{4}\psi^{(0)}\left(  \frac{1}{q-1}-\frac{1}{2}\right)  \left(
H_{\frac{1}{q-1}-\frac{3}{2}}+2\ln\frac{q-1}{\beta^{2}}-2i\pi+\gamma
_{\text{E}}+\ln4\right) \nonumber\\
&  +\frac{1}{2}(\gamma_{\text{E}}-i\pi)\ln\frac{2(q-1)}{\beta^{2}}+\frac{1}%
{4}\ln\frac{q-1}{\beta^{2}}\ln\frac{4(q-1)}{\beta^{2}}\nonumber\\
&  +\frac{1}{4}\psi^{(1)}\left(  \frac{1}{q-1}-\frac{1}{2}\right)  -\frac
{3}{8}\pi^{2}-\frac{1}{2}i\gamma_{\text{E}}\pi+\frac{1}{4}\gamma_{\text{E}%
}^{2}+\frac{1}{4}\ln^{2}2,
\end{align}
and so on.

The above result is consistent with the result obtained from the definition
(\ref{logM}).

\subsection{Normal distribution}

The normal distribution%
\begin{equation}
p\left(  x\right)  =\frac{1}{\sqrt{2\pi}}e^{-x^{2}/2}%
\end{equation}
has $n$-th power moments and does not need renormalization. The $n$-th power
moment of the standard normal distribution can be directly calculated by
definition (\ref{norder}):%
\begin{equation}
m_{n}=\frac{1}{\sqrt{\pi}}2^{\frac{n}{2}-1}\left[  (-1)^{n}+1\right]
\Gamma\left(  \frac{n+1}{2}\right)  .
\end{equation}
By the relation between $n$-th power moments and logarithmic moments, Eq.
(\ref{powerlog}), we obtain the first-order logarithmic moment:%
\begin{equation}
\widetilde{m}_{1}=\left.  \frac{dm_{\sigma-1}}{d\sigma}\right\vert _{\sigma
=1}=-\frac{1}{2}\left(  \ln2+\gamma_{\text{E}}-i\pi\right)
\end{equation}
and the second-order logarithmic moment%
\begin{equation}
\widetilde{m}_{2}=\left.  \frac{d^{2}m_{\sigma-1}}{d\sigma^{2}}\right\vert
_{\sigma=1}=\frac{1}{4}\left(  \ln2+\gamma_{\text{E}}-i\pi\right)  ^{2}%
-\frac{\pi^{2}}{8}.
\end{equation}
The expressions of the third-order and fourth-order logarithmic moments are
too long to be listed here.

The above result is consistent with the result obtained from the definition
(\ref{logM}).

\subsection{Student's t-distribution}

The Student's t-distribution%
\begin{equation}
p\left(  x\right)  =\frac{\left(  \frac{\nu}{\nu+x^{2}}\right)  ^{\frac{\nu
+1}{2}}}{\sqrt{\nu}B\left(  \frac{\nu}{2},\frac{1}{2}\right)  }%
\end{equation}
has the $n$-th power moment, and does not need renormalization. The $n$-th
power moment can be directly calculated by definition (\ref{norder}):%
\begin{equation}
m_{n}=\frac{\left[  (-1)^{n}+1\right]  \nu^{n/2}\Gamma\left(  \frac{n+1}%
{2}\right)  \Gamma\left(  \frac{\nu-n}{2}\right)  }{2\sqrt{\pi}\Gamma\left(
\frac{\nu}{2}\right)  }.
\end{equation}
By the relation between $n$-th power moments and logarithmic moments, Eq.
(\ref{powerlog}), we obtain the first logarithmic moment,%
\begin{equation}
\widetilde{m}_{1}=\left.  \frac{dm_{\sigma-1}}{d\sigma}\right\vert _{\sigma
=1}=\frac{1}{2}\left(  \ln\frac{\nu}{4}+i\pi-H_{\frac{\nu}{2}-1}\right)
\end{equation}
and the second logarithmic moment,%
\begin{align}
\widetilde{m}_{2}  &  =\left.  \frac{d^{2}m_{\sigma-1}}{d\sigma^{2}%
}\right\vert _{\sigma=1}\nonumber\\
&  =\frac{1}{8}\left[  2\ln^{2}\nu+4i\pi\ln\nu+2\psi^{(1)}\left(  \frac{\nu
}{2}\right)  -3\pi^{2}\right. \nonumber\\
&  \left.  +2\left(  \psi^{(0)}\left(  \frac{\nu}{2}\right)  +\gamma
_{\text{E}}+2\ln2\right)  \left(  2\ln2-2\ln\nu+\psi^{(0)}\left(  \frac{\nu
}{2}\right)  -2i\pi+\gamma_{\text{E}}\right)  \right]  .
\end{align}
The expressions of the third-order and fourth-order logarithmic moments are
too long to be listed here.

The above result is consistent with the result obtained from the definition
(\ref{logM}).

\subsection{Laplace distribution}

The probability density function of the Laplace distribution is%
\begin{equation}
p\left(  x\right)  =\frac{\lambda}{2}e^{-\lambda\left\vert x-\mu\right\vert
}=\left\{
\begin{array}
[c]{c}%
\frac{\lambda}{2}e^{-\lambda\left(  x-\mu\right)  },\text{ }x>\mu\\
\frac{\lambda}{2}e^{-\lambda\left(  -x+\mu\right)  },\text{ }x\leq\mu
\end{array}
\right.  .
\end{equation}
The $n$-th power moment of the Laplace distribution by the $n$-th moment
definition (\ref{norder}) is
\begin{equation}
m_{n}=\frac{(-1)^{n}e^{-\lambda\mu}\Gamma(n+1,-\lambda\mu)+e^{\lambda\mu
}\Gamma(n+1,\lambda\mu)}{2\lambda^{n}}.
\end{equation}
By the relation between $n$-th power moments and logarithmic moments, Eq.
(\ref{powerlog}), we obtain the first logarithmic moment,%
\begin{align}
\widetilde{m}_{1}  &  =\left.  \frac{dm_{\sigma-1}}{d\sigma}\right\vert
_{\sigma=1}\nonumber\\
&  =\frac{1}{2}\left[  e^{-\lambda\mu}\Gamma(0,-\lambda\mu)+e^{\lambda\mu
}\Gamma(0,\lambda\mu)\right]  +\ln\mu+i\pi,
\end{align}
and the second logarithmic moment
\begin{align}
\widetilde{m}_{2}  &  =\left.  \frac{d^{2}m_{\sigma-1}}{d\sigma^{2}%
}\right\vert _{\sigma=1}\nonumber\\
&  =e^{-\lambda\mu}G_{2,3}^{3,0}\left(  -\lambda\mu\left\vert
\begin{array}
[c]{c}%
1,1\\
0,0,0
\end{array}
\right.  \right)  +e^{\lambda\mu}G_{2,3}^{3,0}\left(  \lambda\mu\left\vert
\begin{array}
[c]{c}%
1,1\\
0,0,0
\end{array}
\right.  \right) \nonumber\\
&  +e^{\lambda\mu}\Gamma(0,\lambda\mu)\ln\mu+e^{-\lambda\mu}(\ln\mu
+2i\pi)\Gamma(0,-\lambda\mu)+(\ln\mu+i\pi)^{2}-\pi^{2}%
\end{align}
where $G_{p,q}^{m,n}\left(  z\left\vert
\begin{array}
[c]{c}%
a_{1},\cdots,a_{p}\\
b_{1},\cdots,b_{q}%
\end{array}
\right.  \right)  $ is the Meijer $G$ function and $_{p}F_{q}(a;b;z)$ is the
generalized hypergeometric function (\cite{olver2010nist}). The expressions of
the third-order and fourth-order logarithmic moments are too long to be listed here.

The above result is consistent with the result obtained from the definition
(\ref{logM}).

\section{Discussion and conclusion \label{conclusion}}

In this paper, we suggest renormalization schemes for seeking moments for
distributions that have no moments. The key to renormalization is that the
renormalized result must be renormalization scheme-independent. In order to
show that the renormalized moment is scheme-independent, we construct a
variety of renormalization schemes. For each distribution, we use different
renormalization schemes to calculate the same renormalized power moment to
exemplify scheme independence.

The renormalization schemes in this paper can be divided into the following cases.

1) The renormalized moment can be obtained by directly analytically continuing
the moment definition. In this case, the divergence is caused by the form of
the definition. The moment of the distribution is usually defined by
integrals. If it is not integrable, then the distribution has no moment. It is
known that different representations of the same function usually have
different domains. As an analogy, we take the gamma function as an example.
The integral definition of the gamma function is
\begin{equation}
\Gamma\left(  z\right)  =\int_{0}^{\infty}t^{z-1}e^{-t}dt,\text{
\ }\operatorname{Re}z>0.
\end{equation}
In this integral definition, the domain of the gamma function is
$\operatorname{Re}z>0$. If we define the gamma function by the functional
equation,
\begin{equation}
\Gamma\left(  z+1\right)  =z\Gamma\left(  z\right)  ,
\end{equation}
then the gamma function is analytic throughout the whole complex plane except
for isolated singularities. That is to say, under the functional equation
definition, the gamma function is only undefined on some isolated points. The
integral definition and the function equation definition of the gamma function
are equivalent for $\operatorname{Re}z>0$. However, the domains in these two
definitions are different, and the function equation definition gives a larger
domain. That is, through analytical continuation, we can use one definition to
replace another and extend the domain of the gamma function. In addition to
the integral definition and the function equation definition, the gamma
function also has other definitions. These different definitions, which are
analytic continuations to each other, will give the same value in the
overlapping area of their domains. In practice, working out the integral or
the sum in the definition is just to do analytic continuation. The uniqueness
of analytical continuation ensures that the renormalized result is scheme-independent.

2) If analytic continuation of the moment definition still cannot give a
finite moment, which shows that the divergence is caused by true singularities
--- the obstacle of analytic continuation, it requires minimal subtraction to
remove the singularity. The analytical continuation can only remove the
singularities caused by the form of definition, but not the true
singularities. Taking the gamma function as an example, by replacing the
integral definition with the function equation definition, we can extend the
domain of the gamma function to the entire complex plane except for several
isolated singularities. These isolated singularities are true singularities,
which are obstacles to the analytical continuation and cannot be removed by an
analytical continuation. In fact, it is because of these true singularities
that the gamma function, such as integral definition and series definition, is
limited to a certain local area of the complex plane. This means that even if
we extend the domain of moments through analytical continuation, we still
cannot define moments on these true singularities. The classification of
isolated singularities is done through series. Thus we can expose
singularities through series expansion. Furthermore, minimal subtraction is
used to remove these singularities and give the renormalized moments on these
true singularities.

3) A new parameter is introduced, making the moment a function of this
parameter. The divergence of the moment then becomes the singularity of this
function and then is removed by renormalization treatments. In other words,
the integral in the moment definition is divergent. After introducing a new
parameter, the integral only diverges at some special values of this parameter
but is well-defined at other parameter values. This turns the divergence
problem into a problem that can be handled by renormalization. Both the
weighted moment method and the cut-off method belong to this case.

We develop a method to express logarithmic moments by power moments. In this
way, the renormalization treatment of power moments can be applied to
calculating logarithmic moments. This method can also be used to verify the
validity of the renormalization scheme by comparing the logarithmic moment
obtained directly from the definition and from the renormalized power moment.

In addition, technically, the renormalization procedure often relies on the
explicit result of an integral. In different renormalization schemes, we
encounter different integrals. If an integral cannot be worked out in a
renormalization scheme, we can turn to another renormalization scheme. The
more renormalization schemes are, the greater the possibility of successful
renormalization is.

We will discuss the various random processes with renormalized moments in
further consideration.

\acknowledgments

We are very indebted to Dr G. Zeitrauman for his encouragement. This work is supported in part by Special Funds for theoretical physics Research Program of the NSFC under Grant No. 11947124, and NSFC under Grant Nos. 11575125 and 11675119.






\end{CJK*}
\end{document}